\documentclass[dvipsnames]{amsart}

\usepackage{amsthm}
\usepackage{amssymb}
\usepackage{amsmath}
\usepackage{mathtools} 
\usepackage[shortlabels]{enumitem}
\usepackage{tikz-cd}
\usepackage{subfigure}

\usepackage[english]{babel}
\usepackage[shortlabels]{enumitem}
\usepackage{csquotes}

\usepackage[style=alphabetic, backend=biber, sorting=nyt, doi=false, isbn=false,url=false, hyperref,backref,backrefstyle=none, maxnames=20, maxalphanames=20]{biblatex}
\newbibmacro{string+doiurlisbn}[1]{%
  \iffieldundef{doi}{%
    \iffieldundef{url}{%
      \iffieldundef{isbn}{%
        \iffieldundef{issn}{%
          #1%
        }{%
          \href{https://books.google.com/books?vid=ISSN\thefield{issn}}{#1}%
        }%
      }{%
        \href{https://books.google.com/books?vid=ISBN\thefield{isbn}}{#1}%
      }%
    }{%
      \href{\thefield{url}}{#1}%
    }%
  }{%
    \href{https://dx.doi.org/\thefield{doi}}{#1}%
  }%
}
\DeclareFieldFormat{title}{\usebibmacro{string+doiurlisbn}{\mkbibemph{#1}}}
\DeclareFieldFormat[article,incollection,inproceedings,unpublished,misc,book]{title}%
    {\usebibmacro{string+doiurlisbn}{\mkbibquote{#1}}}
\DefineBibliographyStrings{english}{%
  backrefpage = {$\uparrow$},%
  backrefpages = {$\uparrow$}%
}

\usepackage{url}
\usepackage[hyperindex,breaklinks]{hyperref}

\theoremstyle{plain}
\newtheorem{theorem}{Theorem}[section]
\newtheorem{proposition}[theorem]{Proposition}
\newtheorem{lemma}[theorem]{Lemma}
\newtheorem{corollary}[theorem]{Corollary}

\theoremstyle{definition}
\newtheorem{definition}[theorem]{Definition}
\newtheorem{remark}[theorem]{Remark}
\newtheorem{example}[theorem]{Example}
\newtheorem{algorithm}[theorem]{Algorithm}
\numberwithin{equation}{section}

\def\Q{\mathbb{Q}}
\def\Z{\mathbb{Z}}
\def\F{\mathbb{F}}

\def\C{\mathbb{C}}

\newcommand{\ZFV}{\mathbb{Z}[\pi,\overline\pi]}

\DeclareMathOperator{\len}{\mathrm{length}}
\DeclareMathOperator{\ICM}{ICM}
\DeclareMathOperator{\Pic}{Pic}
\DeclareMathOperator{\Hom}{Hom}
\DeclareMathOperator{\End}{End}
\DeclareMathOperator{\Tr}{Tr}

\DeclareMathOperator{\dom}{dom}
\DeclareMathOperator{\cod}{cod}
\DeclareMathOperator{\image}{image}

\newcommand{\fra}{\mathfrak{a}}
\newcommand{\frb}{\mathfrak{b}}
\newcommand{\frc}{\mathfrak{c}}
\newcommand{\frf}{\mathfrak{f}}

\newcommand{\fri}{\mathfrak{i}}
\newcommand{\frj}{\mathfrak{j}}
\newcommand{\frk}{\mathfrak{k}}
\newcommand{\frl}{\mathfrak{l}}
\newcommand{\frm}{\mathfrak{m}}
\newcommand{\frp}{\mathfrak{p}}
\renewcommand{\frq}{\mathfrak{q}}
\newcommand{\frw}{\mathfrak{w}}
\newcommand{\cE}{{\mathcal E}}
\newcommand{\cF}{{\mathcal F}}
\newcommand{\cG}{{\mathcal G}}
\newcommand{\cI}{{\mathcal I}}
\newcommand{\cO}{{\mathcal O}}
\newcommand{\cP}{{\mathcal P}}
\newcommand{\cR}{{\mathcal R}}
\newcommand{\cS}{{\mathcal S}}
\newcommand{\cU}{{\mathcal U}}
\newcommand{\cV}{{\mathcal V}}
\newcommand{\cW}{{\mathcal W}}
\newcommand{\cL}{{\mathcal L}}

\newcommand{\ucirc}{\underline{\circ}}
\newcommand{\usq}{\underline{\bullet}}

\newcommand{\avlink}[1]{\href{http://www.lmfdb.org/Variety/Abelian/Fq/#1}{\textsf{#1}}}

\title{Ordinary abelian varieties: isogeny graphs and polarizations}
\date{\today}

\author[E. Costa]{Edgar Costa}
\address{Edgar Costa, Department of Mathematics, Massachusetts Institute of Technology, Cambridge, MA 02139-4307, USA}
\email{edgarc@mit.edu}
\urladdr{\url{https://edgarcosta.org}}

\author[T. Dupuy]{Taylor Dupuy}
\address{Taylor Dupuy, Department of Mathematics and Statistics,
	Innovation Hall E220,
	82 University Place,
	Burlington, VT, 05405, USA}
\email{taylor.dupuy@gmail.com}
\urladdr{\url{http://tdupu.github.io}}

\author[S. Marseglia]{Stefano Marseglia}
\address{Stefano Marseglia, Laboratoire Jean Alexandre Dieudonné, Université Côte d'Azur, 06108 Nice Cedex 2, France}
\email{stefano.marseglia@univ-cotedazur.fr}
\urladdr{\url{https://stmar89.github.io}}

\author[D. Roe]{David Roe}
\address{David Roe, Department of Mathematics, Massachusetts Institute of Technology, Cambridge, MA 02139-4307, USA}
\email{roed@mit.edu}
\urladdr{\url{https://math.mit.edu/~roed}}

\author[C. Vincent]{Christelle Vincent}
\address{Christelle Vincent, Department of Mathematics and Statistics,
	Innovation Hall E220,
	82 University Place,
	Burlington, VT, 05405, USA}
\email{christelle.vincent@uvm.edu}
\urladdr{\url{https://www.uvm.edu/~cvincen1}}

\begin{document}

\begin{abstract}
Given an integer $D$ and an ordinary isogeny class of abelian varieties defined over a finite field $\F_q$ with commutative $\F_q$-endomorphism algebra, we provide algorithms for computing all isogenies of degree dividing $D$ and polarizations of degree dividing $D$.  We discuss phenomena that arise for higher dimension abelian varieties but not elliptic curves, bounds on the diameter of the graph of minimal isogenies, and decompositions of isogeny graphs into orbits for the Picard group of the Frobenius order.
\end{abstract}

\maketitle

\section{Introduction}
In his PhD thesis, Kohel \cite{Kohelthesis} pioneered the study of the structure of isogeny graphs of elliptic curves defined over finite fields $\F_q$.
For ordinary elliptic curves, these graphs have a striking form, called ``isogeny volcanoes'' \cite{FouquetMorain}.
We refer the interested reader to \cite{Sutherland_volcanoes} for details and a more complete history.
A natural generalization is to consider the case of higher-dimensional abelian varieties, and the first results in this direction were obtained for the case of absolutely simple ordinary abelian varieties in \cite{BJW17,Martindalethesis,IonicaThome}.
More recently, \cite{ArpinMarsegliaSpringer_arXiv} obtain results for unpolarized abelian varieties with commutative $\F_q$-endomorphism algebra.

Here we continue with the study of the isogeny graphs of ordinary abelian varieties, and apply this work to the computation of the polarizations of these varieties. A goal of this work, along with \cite{lucant}, is for the eventual inclusion of isomorphism classes of abelian varieties defined over finite fields in the L-functions and modular forms database
\cite{lmfdb}, along with their isogeny graphs.

More precisely, let $q$ be a power of a prime and $h(x)$ be a $q$-Weil polynomial.
By Honda-Tate theory  \cite{Tate66,Honda68}, such a polynomial determines an $\F_q$-isogeny class $\cI_h$ of ordinary abelian varieties defined over $\F_q$.
We further assume $h(x)$ has no repeated complex roots, or equivalently that each abelian variety in $\cI_h$ has commutative $\F_q$-endomorphism algebra.
In this setting, we study the \emph{$D$-isogeny graph}  of $\cI_h$ for $D > 1$, whose vertices are the $\F_q$-isomorphism classes of abelian varieties in $\cI_h$ and whose edges are \emph{equivalence classes} of $\F_q$-isogenies of degree dividing $D$; see Definitions~\ref{def:isog_equiv} and \ref{def:graphs} for more details.
Throughout, all morphisms, isomorphism classes, and isogenies are defined over $\F_q$.

In Section~\ref{sec:icm}, we use Deligne's functor to connect abelian varieties with ideals in the étale algebra $K=\Q[x]/(h(x))$.
Denoting by $\pi$ the class of $x$ in $K$ and by $\overline \pi$ its complex conjugate,
the order $\ZFV \subset K$ is the \emph{Frobenius order} of the isogeny class $\cI_h$.
Its set of fractional ideals is a commutative monoid under ideal multiplication, and isomorphism classes of fractional $\ZFV$-ideals form the \emph{ideal class monoid} $\ICM(\ZFV)$ \cite[Corollary 3.4]{MarICM18}.
The orbits of the action of $\Pic(\ZFV)$ on $\ICM(\ZFV)$ by ideal multiplication are called \emph{weak equivalence classes}.

In Section~\ref{sec:av} we start studying isogenies in $\cI_h$.
The main results of this section show that minimal isogenies necessarily have prime-power degree (Proposition \ref{prop:isog}) and that one of the endomorphism ring of the domain or codomain is contained in the endomorphism ring of the other (Theorem \ref{thm:prop_min_isog}), but there are examples where this inclusion is not maximal (Example \ref{ex:notmin}).
This contrasts with the case of $\ell$-isogenies of elliptic curves for $\ell$ prime, where the inclusion is always maximal.

We continue in Section~\ref{sec:graph} by formally defining the $D$-isogeny graphs $\cG^D$.
Proposition \ref{prop:missing_edges} gives a necessary condition for the existence of minimal isogenies of given prime-power degree; Theorem \ref{thm:conn} provides a sufficient condition for as well as diameter bounds.
We then give in Section~\ref{sec:algorithms_GD} an algorithm to compute $\cG^D$ by leveraging the action of $\Pic(\ZFV)$ to compute minimal edges for each equivalence class. Composing these minimal isogenies then yields $\cG^D$; see Algorithm \ref{alg:compositions}.

In Section~\ref{sec:alg_storing}, we give an algorithm to efficiently store $\cG^D$.
Lemma \ref{lem:GST_free} shows that a quotient of $\Pic(\ZFV)$ acts freely on edges between weak equivalence classes; Algorithm \ref{alg:GST_orbits} computes orbit representatives.
These results allow storage of only one edge from each orbit, which can both compactify outputs and speed up computations.
Finally, in Section~\ref{sec:pols}, we show how to compute polarizations of degree dividing $D$ in $\cI_h$ using $\cG^D$ as input; see Algorithm \ref{alg:pols_from_dualiso}.

All of the code accompanying this paper is available at \url{https://github.com/stmar89/AbVarFqIsogenies}. We conclude with an example of Algorithm \ref{alg:compositions}.

\begin{example}\label{example:4.3.c-ab-af-ai}
Let $D=2$ and consider the ordinary isogeny class of abelian $4$-folds over $\F_3$ with LMFDB label \avlink{4.3.c\_ab\_af\_ai} and Weil polynomial $h(x)=x^8 + 2x^7- x^6 - 5x^5 - 8x^4 - 15x^3 - 9x^2 + 54x + 81.$
This isogeny class has $14$ isomorphism classes; hence its $2$-isogeny graph (Figure \ref{figure:4.3.c-ab-af-ai}) has $14$ vertices, which are here colored by endomorphism ring (darker colors indicating a larger ring).
\begin{figure}[!ht]
    \begin{center}
        \includegraphics[scale=0.33]{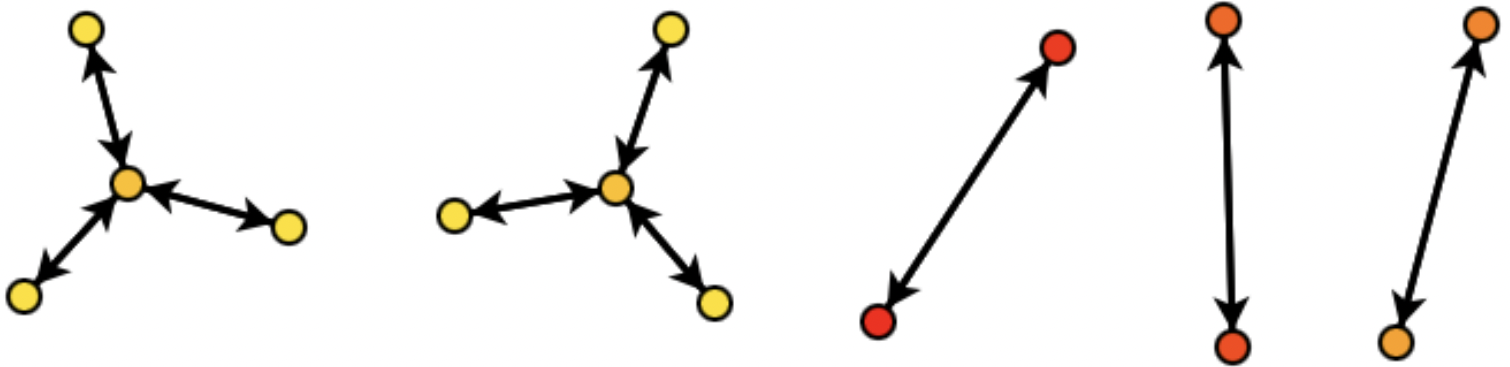}
    \end{center}
    \caption{The 2-isogeny graph for \avlink{4.3.c\_ab\_af\_ai}.}\label{figure:4.3.c-ab-af-ai}
\end{figure}

The base change of this isogeny class to $\F_9$ has Weil polynomial
\[ h_2(x)=x^8 - 6x^7 + 5x^6 + 33x^5 - 122x^4 + 297x^3 + 405x^2 - 4374x + 6561, \]
which is not in the LMFDB, but can be computed as in \cite[\S 3.4]{LMFDB_paper}.
This isogeny class has $1749$ isomorphism classes. Its $2$-isogeny graph has $35$ components: two with $42$ vertices, nine with $17$, twelve with $94$, and twelve with $32$ vertices.

Interestingly, in both cases all components of the same size are isomorphic as directed graphs, but the vertices need not have identical coloring. This is illustrated by one of the component with 2 vertices in Figure \ref{figure:4.3.c-ab-af-ai}, and occurs again over $\F_9$: for example, six of the components with 17 vertices are isomorphic as colored directed graphs, but the last three such components are all pairwise colored differently; two of the four different colorings are shown in Figure \ref{figure:4.3.c-b-af-ai--2--17--1.png}.

\begin{figure}[!ht]
\centering
\begin{subfigure}
  \centering
  \includegraphics[scale=0.23]{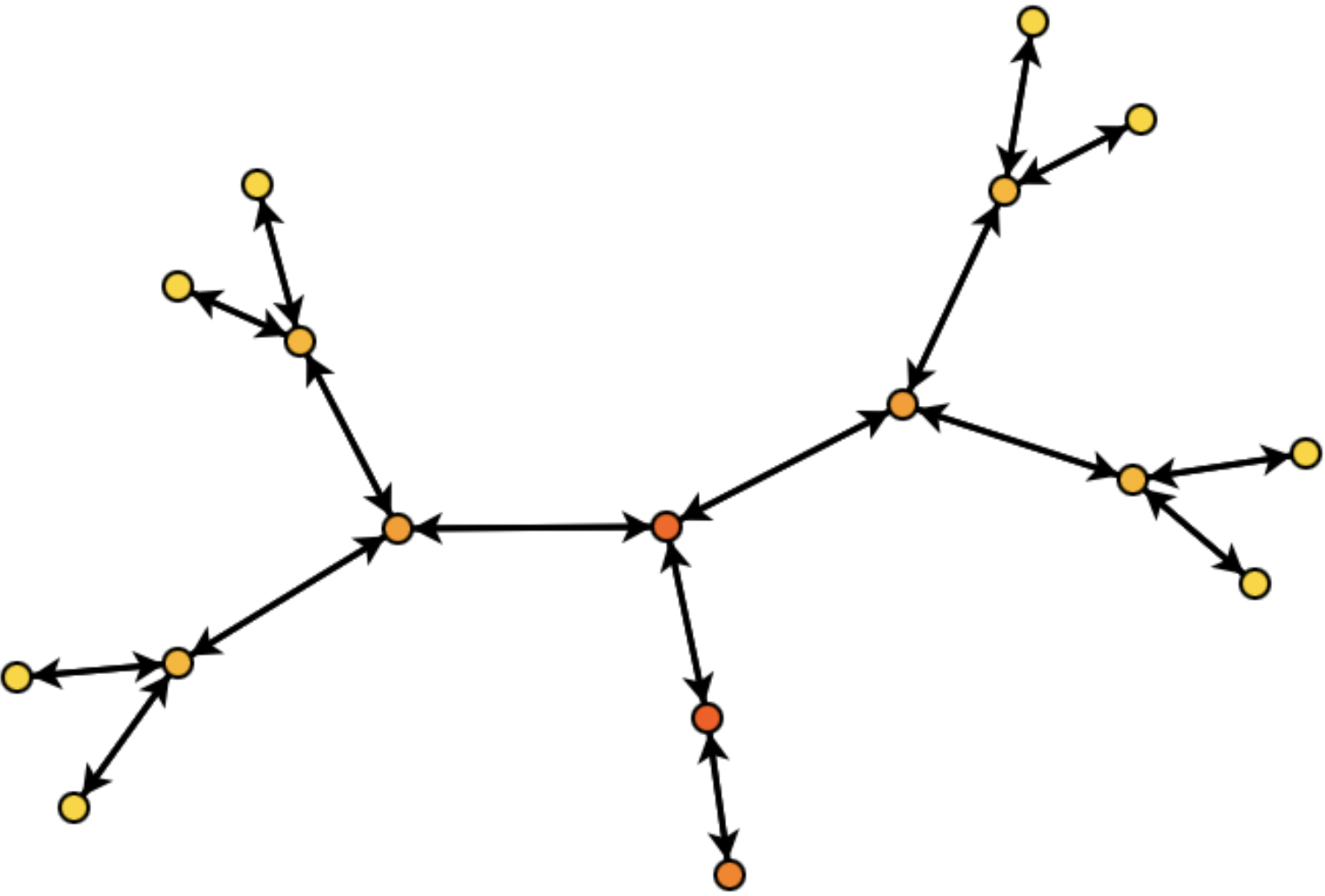}
\end{subfigure}%
\begin{subfigure}
  \centering
\includegraphics[scale=0.23]{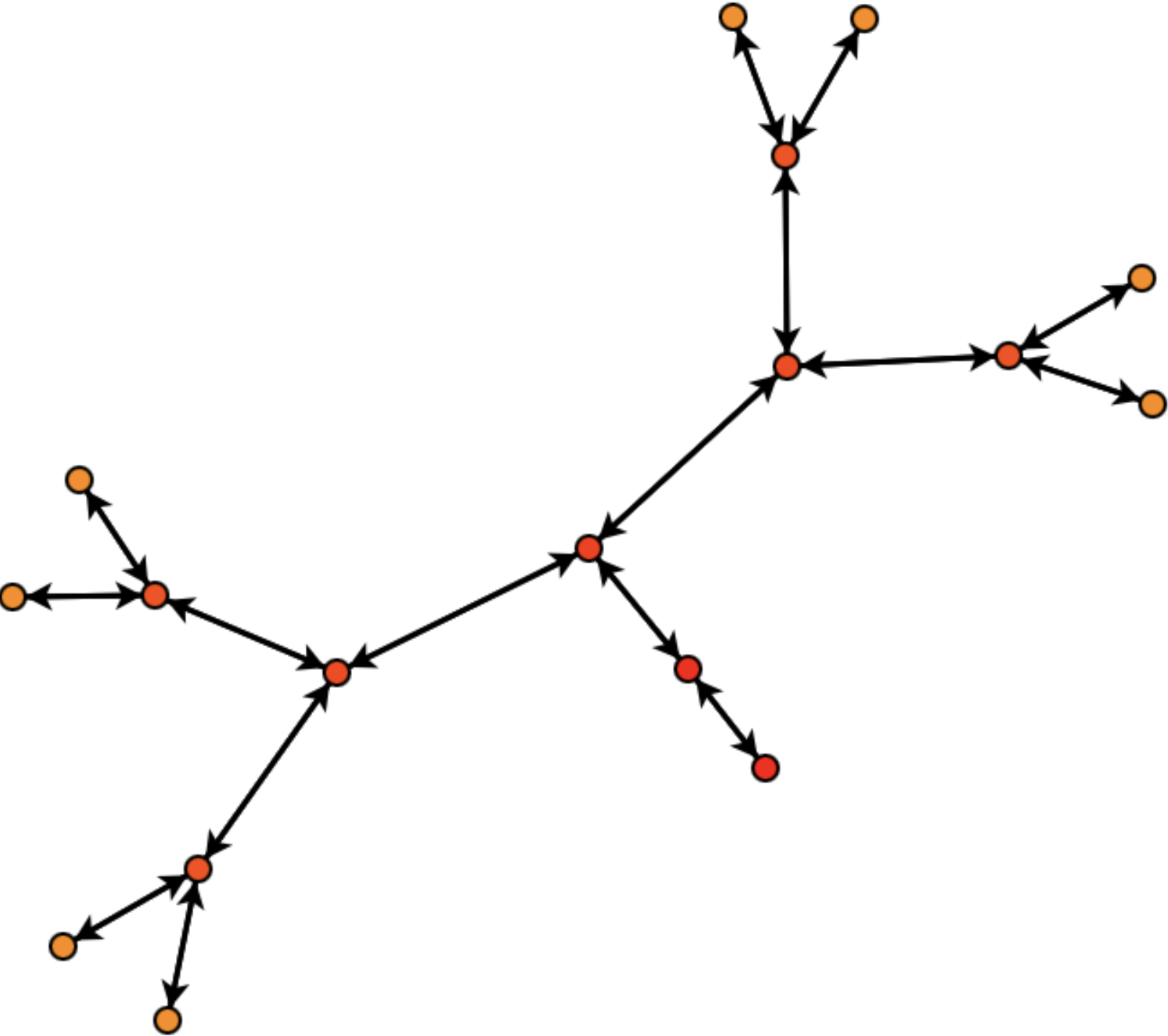}
\end{subfigure}
\caption{Two components of the $2$-isogeny graph of the base change of \avlink{4.3.c\_ab\_af\_ai} to $\F_9$, each with 17 vertices and isomorphic as directed graphs, but whose underlying vertices have different endomorphism rings.}\label{figure:4.3.c-b-af-ai--2--17--1.png}
\end{figure}

Lest the reader be left with the impression that these isogeny graphs are classical isogeny volcanoes, we show a component with 94 vertices in Figure~\ref{figure:4.3.c-b-af-ai--2--94.png}. 
The vertices of this graph have six different endomorphism rings, which in this case correspond to different levels in the graph. The yellow vertices are abelian varieties whose endomorphism ring is the Frobenius order, which is necessarily the smallest order, and the largest endomorphism ring contains the Frobenius order with index $32$ and is contained in the maximal order of $K$ with index $625=5^4$. This component does not contain any horizontal isogenies.
\begin{figure}[!ht]
    \begin{center}
        \includegraphics[scale=0.25]{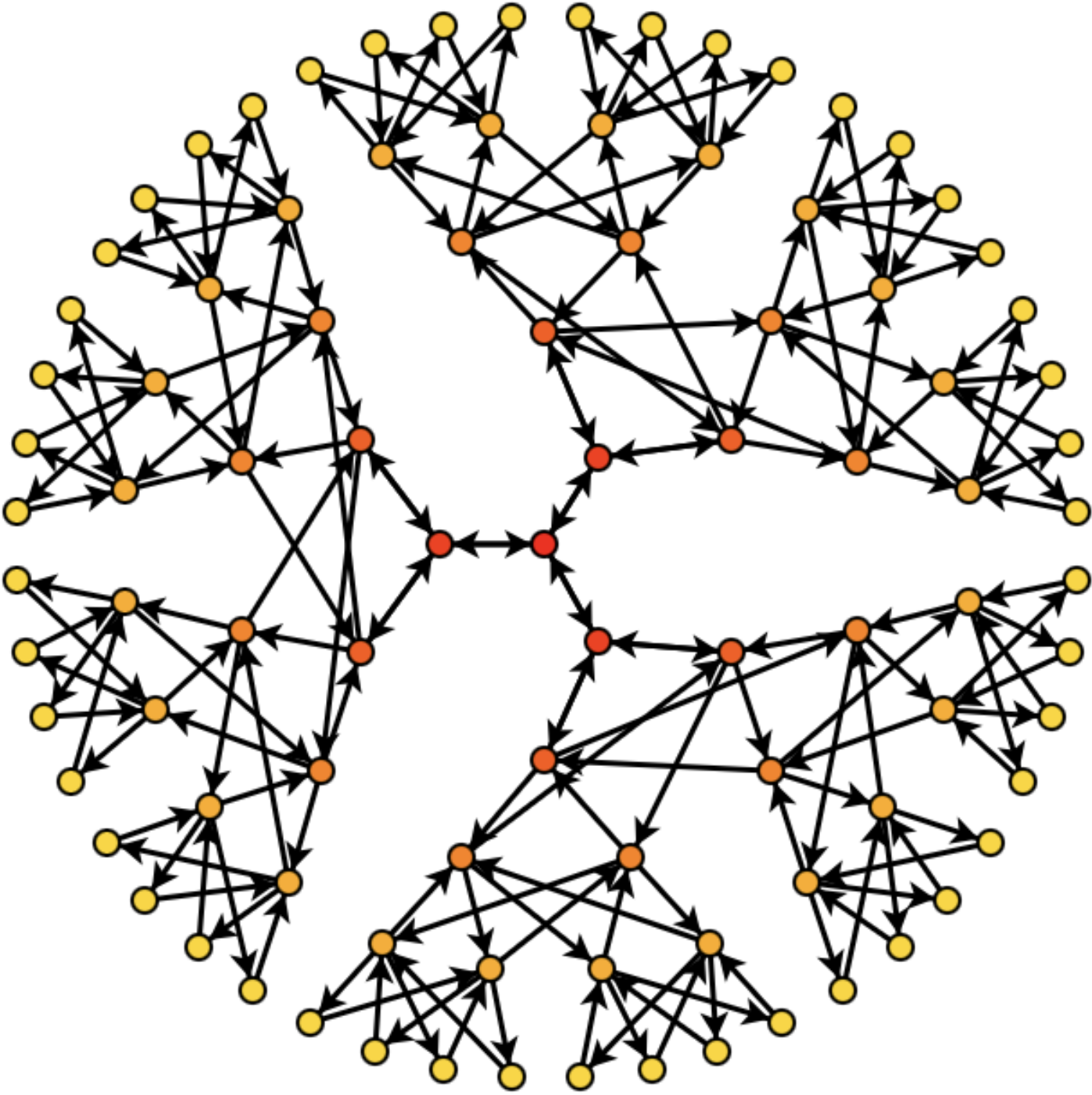}
    \end{center}
    \caption{A component of the 2-isogeny graph for the base change of  \avlink{4.3.c\_ab\_af\_ai} to $\F_9$ with 94 vertices.}
    \label{figure:4.3.c-b-af-ai--2--94.png}
\end{figure}
\end{example}

\section{The ideal class monoid}\label{sec:icm}

In \cite{Del69}, Deligne constructed an equivalence between the category of abelian varieties in an $\F_q$-isogeny class $\cI_h$, together with their $\F_q$-homomorphisms, and the category of fractional $\ZFV$-ideals, together with $\ZFV$-linear morphisms.
We recall the construction, as we need the notation below.

To define the functor, we use that an ordinary abelian variety $A$ defined over a finite field $\F_q$ of characteristic $p$, together with the action of its Frobenius endomorphism, admits a lifting $\tilde A$  to a finite extension of $\Q_p$, called the canonical lift of $A$.
Deligne's functor is then defined by the association $A\mapsto T(A)=H_1(\tilde A\otimes_\varepsilon \C,\Z)$, where $\varepsilon \colon \overline{\Q}_p\to \C$ is any fixed isomorphism.
The $\Z$-module $T(A)$ is finitely generated and free of rank $2\dim(A)$, comes equipped with natural $\Z$-linear morphisms $\pi$ and $q/\pi$ which correspond to the Frobenius and Verschiebung endomorphisms of $A$, and can be considered as a global analogue of the Tate module of $A$.
In particular, this module allows us to study isogenies from and to $A$, including those where $p$ divides the degree.
The restriction $\cF$ of this functor to $\cI_h$ gives a bijection between isomorphism classes and $\ICM(\ZFV)$, allowing us to represent varieties and isogenies via fractional $\ZFV$-ideals \cite{MarAbVar18}.

Let $R$ be an order in an \'etale algebra $K$. 
For each pair of fractional $R$-ideals $I,J$, we define $(I:J)=\{a\in K: aJ\subseteq I\}$, the \emph{colon ideal} of $I$ and $J$.
Note that $(I:I)$ is an overorder of $R$, called the \emph{multiplicator ring} of $I$.
As in the case of the Frobenius order $\ZFV$, ideal multiplication induces the structure of a commutative monoid on the set of fractional $R$-ideals modulo $R$-linear isomorphisms, and we again call it the \emph{ideal class monoid} of $R$ and denote it by $\ICM(R)$.
For each overorder $S$ of $R$, the set of invertible fractional $S$-ideals forms an abelian group, called the \emph{Picard group} $\Pic(S)$ of $S$, and we have a natural inclusion $\Pic(S)\subseteq\ICM(R)$.
This inclusion is an equality if and only if $S=R=\cO_K$ is the maximal order of $K$; in this case $\ICM(\cO_K)$ is the direct product of class groups.

If $\frm$ is a maximal $R$-ideal, we denote by $R_\frm$ the completion of $R$ at $\frm$ and, for a fractional $R$-ideal $I$, we write $I_{\frm} = I \otimes_R R_{\frm}$.
Two fractional $R$-ideals $I,J$ are \emph{weakly equivalent} if $I_\frm\simeq J_\frm$ for all maximal $\frm$. 
$I$ and $J$ are weakly equivalent precisely when $(I:I)=(J:J)$ and $(I:J)$ is invertible in $(I:I)$ by \cite[Corollary~4.5]{MarICM18}.
Again, as in the case of $\ZFV$, ideal multiplication induces a commutative monoid structure on the set $\cW(R)$ of weak equivalence classes of fractional $R$-ideals.
The monoids $\ICM(R)$ and $\cW(R)$ can be partitioned by
\[ \ICM(R) = \bigsqcup_S \ICM_S(R) \qquad\text{and}\qquad \cW(R)=\bigsqcup_S\cW_S(R), \]
where $S$ runs over the overorders of $R$, and $\ICM_S(R)$ (resp.~$\cW_S(R)$) contains the ideal (resp.~weak equivalence) classes with multiplicator ring $S$.

\begin{proposition}[{\cite[Theorem~4.6]{MarICM18}}]\label{prop:pic_weak}
    The group $\Pic(S)$ acts freely on $\ICM_S(R)$, and the quotient space of this action is precisely $\cW_S(R)$.
\end{proposition}

The following statement is probably well known to the expert, but we could not locate it in the literature.
We provide a simple direct proof.
\begin{lemma}\label{lem:pic_ext}
    Let $S$ be an overorder of $R$, and $J$ be an invertible fractional $S$-ideal.
    Then there exists an invertible $R$-ideal $I$ such that $IS=J$.
    In particular, the extension map $I\mapsto IS$ 
    induces a surjective group homomorphism
    \[ e_S \colon \Pic(R) \to \Pic(S). \]
\end{lemma}
\begin{proof}
Let $y \in K^\times$ be such that $J'=yJ\subseteq R$, and let $\frp_1,\ldots,\frp_n$ be the support of the finite $R$-module $S/J'$.
For each maximal ideal $\frp$ of $R$, the localization $S_\frp$ is a semilocal ring.
This implies that for each $\frp_i$ in the support there exists an element $x_i\in J'$ such that $J'_{\frp_i}=(x_iS)_{\frp_i}$, while for every other maximal ideal $\frq$ of $R$ we have that $J'_\frq=S_\frq$.
Now, for each $i$, let $k_i$ be an integer satisfying $\frp_i^{k_i}\subseteq (x_iR)_{\frp_i}$.
Define
\[ I' = \sum_{i=1}^n \left( x_iR + \frp_i^{k_i} \right)\prod_{j\neq i} \frp_j^{k_j}. \]
A direct verification shows that $I'_{\frp_i}=(x_iR)_{\frp_i}$ for all $i$ and $I'_{\frq}=R_{\frq}$ for the other maximal ideals $\frq$ of $R$.
In particular, $I'$ is invertible, since it is locally principal, and furthermore $I'S = J'$. Therefore for $I=y^{-1}I'$ we have $IS=J$ as required.
\end{proof}

Finally, for a fractional $R$-ideal $I$, we denote by $I^\dagger$ its \emph{trace dual ideal}
\[ I^\dagger = \{ x\in K : \Tr_{K/\Q}(xI)\subseteq \Z \}. \]

\begin{lemma} \label{lem:dagger} Let $x \in K^\times$, $L$ be an invertible $R$-ideal, and $I$ be a fractional $R$-ideal.
\begin{enumerate}
\item We have $(xI)^\dagger = \frac{1}{x}I^\dagger$ and $(IL)^\dagger = I^\dagger L^{-1}$.
\item The map $I\mapsto I^\dagger$ induces well defined maps $\ICM(R)\to \ICM(R)$ and $\cW(R)\to \cW(R)$.
\end {enumerate}
\end{lemma}
\begin{proof}
    The first statement follows by direct computation and localization, and the second statement from the first.
\end{proof}

\section{Isogenies}\label{sec:av}

We now discuss isogenies between isomorphism classes and set $R = \ZFV$ for the remainder of the paper.
Because we work up to isomorphism, we impose the following equivalence relation:

\begin{definition} \label{def:isog_equiv}
    Two isogenies $\varphi_1 \colon A_1\to B_1$ and $\varphi_2 \colon A_2\to B_2$ are said to be \emph{equivalent}, denoted $\varphi_1 \approx \varphi_2,$ if and only if there exist isomorphisms $i_1\colon A_1\to A_2$ and $i_2\colon B_1\to B_2$ such that $i_2\circ \varphi_1 = \varphi_2\circ i_1$.
\end{definition}

For the purposes of computing the isogeny graphs we seek, we will need the following two notions:

\begin{definition}
    An isogeny $\varphi \colon A\to B$ of degree larger than $1$ is called \emph{minimal} if, given a decomposition $\varphi = \psi_1\circ \psi_2$, we have that either $\psi_1$ or $\psi_2$ is an isomorphism.
\end{definition}

As well as:

\begin{definition}
The \emph{length} of an isogeny $\varphi$ is the smallest positive integer $r$ such that there exists a decomposition $\varphi = \varphi_1 \circ \dots \circ \varphi_r$ of $\varphi$ into minimal isogenies $\varphi_i$.
\end{definition}

Both the length and degree are constant in an equivalence class of isogenies.

In Proposition \ref{prop:isog}, we obtain some properties of the Deligne functor $\cF$ which we will need to study isogenies between abelian varieties in $\cI_h$.

\begin{proposition}\label{prop:isog}
    Let $A,B\in \cI_h$, and $I=\cF(A)$ and $J=\cF(B)$.
    Then $\cF$ sends an isogeny $\varphi\colon A \to B$ to an element $a_\varphi\in (J:I) \cap K^\times$.
    Moreover,
    we have:
    \begin{enumerate}[(i)]
       \item\label{prop:isog:deg}
       $\deg(\varphi) = [ J : a_\varphi I ]$.
       \item\label{prop:isog:incl}
       Given a positive integer $D$, if $\deg \varphi$ divides $D$ then $DJ \subseteq a_\varphi I$.
       \item\label{prop:isog:min}
       The isogeny $\varphi$ is minimal if and only if the inclusion $a_\varphi I \subsetneq J$ is maximal, that is, the quotient $J/a_\varphi I$ is a simple $R$-module.
       \item\label{prop:isog:length}
       The length of the isogeny $\varphi$ equals the length of $J/a_\varphi I$ as an $R$-module.
    \end{enumerate}
\end{proposition}

\begin{proof}
    We use the notation introduced in Section \ref{sec:icm} recalling the construction of Deligne's functor.
    Let $\tilde\varphi \colon \tilde A\to \tilde B$ be the canonical lift of $\varphi$, then $\deg(\varphi) = \deg(\tilde\varphi) = \deg(\tilde \varphi \otimes_\varepsilon \C)$.
   Since $\tilde{\varphi}\otimes \C$ is separable and taking integral homology, we get
    \[ \deg(\varphi) = \left\vert \frac{T(B)}{T(\tilde\varphi)T(A)} \right\vert = [ J:a_\varphi I]. \]
    This concludes the proof of \ref{prop:isog:deg}.
    Part~\ref{prop:isog:incl} is then an immediate consequence of~\ref{prop:isog:deg}, since a finite abelian group is annihilated by its order.

    We now prove~\ref{prop:isog:min}.
    Suppose that $J/a_{\varphi} I$ is not simple, so there exists an $R$-module $M$ such that $a_\varphi I \subsetneq M \subsetneq J$.
    Using Deligne's functor, this is equivalent to the existence of two isogenies $\psi_1\colon A\to C$ and $\psi_2\colon C\to B$ such that $\varphi=\psi_2\circ\psi_1$, representing the inclusions $a_\varphi I \subsetneq M$ and $1 \cdot M\subsetneq J$.
    Since the inclusions are strict, neither $\psi_1$ nor $\psi_2$ is an isomorphism by~\ref{prop:isog:deg}.

    We finally consider~\ref{prop:isog:length}.
    The length of $J/a_\varphi I$ as an $R$-module is the length of a chain of inclusions of $R$-modules
    \[ a_\varphi I = M_n \subsetneq M_{n-1} \subsetneq \ldots \subsetneq M_1 = J \]
    in which every quotient $M_i/M_{i+1}$ is a simple $R$-module.
    The conclusion then follows by functoriality and~\ref{prop:isog:min}.
\end{proof}

Before we state our next result we need the following:

\begin{definition}\label{df:horiz_asc_desc}
    Let $\varphi \colon A\to B$ be an isogeny in $\cI_h$.
    We say that $\varphi$ is \emph{horizontal} if $\End(A)=\End(B)$; $\varphi$ is \emph{ascending} if $\End(A)\subsetneq \End(B)$ and \emph{descending} if $\End(B)\subsetneq \End(A)$; and finally $\varphi$ is \emph{vertical} if it is ascending or descending.
\end{definition}

\begin{theorem}\label{thm:prop_min_isog}
    Let $\varphi \colon A\to B$ be a minimal isogeny in $\cI_h$, $I=\cF(A)$, $J=\cF(B)$, $a=\cF(\varphi)$, and $\frm$ be the unique maximal ideal of $R$ such that $J/aI \simeq R/\frm$.
    Then:
    \begin{enumerate}[(i)]
        \item \label{thm:prop_min_isog:compareEnds} $\varphi$ is horizontal or vertical, that is, $\End(A) \subseteq \End(B)$ or $\End(B)\subseteq \End(A)$.
        \item \label{thm:prop_min_isog:m_reg} If $\frm$ is regular then $I$ and $J$ are weakly equivalent. In particular, $\varphi$ is horizontal.
        \item \label{thm:prop_min_isog:m_sing} If $\varphi$ is vertical then $\frm$ is singular.
    \item \label{thm:prop_min_isog:vert_index} Let $T$ be the larger of the two orders $\End(A)$ and $\End(B)$, and $S$ the smaller of the two.
     Then $\vert T/S\vert = \vert R/\frm\vert^{\len_{R}(T/S)}$.
        In particular, $\varphi$ is vertical if and only if $[R:\frm]$ divides the index $[T:S]$.
    \end{enumerate}
\end{theorem}
\begin{proof}
    Using Deligne's functor we have $\End(A)=(I:I)$ and $\End(B)=(J:J)$, and since $(aI:aI)=(I:I)$ in what follows we replace $aI$ with $I$ for simplicity.
    Furthermore, since $\varphi$ is minimal the finite $R$-module $J/I$ is simple by Proposition~\ref{prop:isog}.\ref{prop:isog:min}, and there is a unique maximal $R$-ideal $\frm$ such that $J/I\simeq R/\frm$.

    Now suppose that $\End(B) \not\subseteq \End(A)$, so there is an element $\beta \in (J:J)\setminus (I:I)$.
    Therefore we have $\beta J\subseteq J$ but $\beta I\not\subseteq I$, giving the inclusions
    $ I \subsetneq I+\beta I \subseteq J $.
    Since $J/I$ is a simple $R$-module, we get $I+\beta I = J$.
    Hence, for every $b\in J$ there are $a_0,a_1\in I$ such that $b=a_0+\beta a_1$, and thus for every $\alpha\in (I:I)$, we have
    \[ \alpha b = \alpha a_0 + \beta \alpha a_1 \in I+\beta I = J.\]
    This shows that $\alpha \in (J:J)$, and so $(I:I)\subseteq (J:J)$, thus proving part~\ref{thm:prop_min_isog:compareEnds}.

    Since $J/I\simeq R/\frm$, we get that $I$ and $J$ are locally equal at every maximal ideal other than $\frm$.
    If $\frm$ is regular, then $R_\frm = \cO_{K,\frm}$, which is a product of discrete valuation rings.
    It follows that $I_\frm$ and $J_\frm$ are isomorphic, so $I$ and $J$ are weakly equivalent.
    In particular, $(I:I)=(J:J)$, that is, $\varphi$ is a horizontal isogeny.
    This completes the proof of part \ref{thm:prop_min_isog:m_reg}, and part~\ref{thm:prop_min_isog:m_sing} is immediate.

    We finally prove part~\ref{thm:prop_min_isog:vert_index}.
    As observed above $I$ and $J$ are locally equal at every maximal ideal other than $\frm$.
    Since taking colon ideals commutes with localization, the same is true for their multiplicator rings.
    Hence, $T/S$ is either trivial or a finite $R$-module whose support contains only $\frm$.
    In both cases, we have the equality $\vert T/S\vert = \vert R/\frm\vert^{\len_{R}(T/S)}$.
    The former case happens precisely when $\varphi$ is horizontal, thus completing the proof.
\end{proof}

In the $\ell$-isogeny graph of ordinary elliptic curves, if the isogeny $\varphi \colon A\to B$ is vertical then the index of one of the endomorphism rings in the other is at most $\ell$.
This implies that there is no order strictly between $\End(A)$ and $\End(B)$ if $\ell$ is prime.
The following shows that this results fails in general for minimal isogenies.

\begin{example}\label{ex:notmin}
    Let $\cI_h$ be the isogeny class of abelian surfaces defined over $\F_5$ with LMFDB label \avlink{2.5.a\_g} given by the Weil polynomial $h(x)=x^4+6x^2+25 = (x^2 - 2x + 5)(x^2 + 2x + 5)$.
    Let $R$ be the Frobenius order of $\cI_h$, which has index $[\cO_K:R]=64$ in the maximal order $\cO_K$ of the endomorphism algebra $K$.
    Then $R$ has a unique singular maximal ideal $\frm$ with residue field $\F_2$, as well as a unique overorder $S$ satisfying $[\cO_K:S]=8$.
    Furthermore, the conductor $\frf_S=(S:\cO_K)$ of $S$ in $\cO_K$ is the unique singular maximal ideal of $S$, we have an isomorphism $S/\frf_S\simeq R/\frm$, and $(\frf_S:\frf_S) = \cO_K$.

    Let $\varphi \colon A\to B$ be the descending isogeny that is sent to the inclusion $1\cdot \mathfrak{f}_S \subset S$ by Deligne's functor; this isogeny is minimal of degree $2$ by Proposition~\ref{prop:isog} parts \ref{prop:isog:deg} and \ref{prop:isog:min} since $S/\mathfrak{f}_S$ is simple of size 2.
    Relating this to Theorem~\ref{thm:prop_min_isog}.\ref{thm:prop_min_isog:vert_index}, we observe that $\len_R(\cO_K/S) = 3 = 2\dim(A)-1$, and indeed $[\End(A):\End(B)] = [\cO_K:S] = \deg(\varphi)^{3}$.
    This agrees with \cite[Lemma 2.1]{Gaetan} where isogenies with kernel $(\Z/\ell\Z)^g$ are considered instead.

    Strikingly, there are several orders $T$ such that $S\subsetneq T \subsetneq \cO_K$, despite the fact that $\varphi$ is minimal.
    Indeed, $S$ has $4$ overorders with index $4$ in $\cO_K$ and $3$ overorders with index $2$ in $\cO_K$.
\end{example}

\begin{proposition}\label{prop:same_edge}
	Let $\varphi_1\colon A_1\to B_1$ and $\varphi_2\colon A_2\to B_2$ be two isogenies in $\cI_h$.
	For $i=1,2$, let $I_i = \cF(A_i)$, $J_i = \cF(B_i)$, $x_i = \cF(\varphi_i)$ and $M_i= x_iI_i \subseteq J_i$.
    The following are equivalent:
    \begin{enumerate}[(i)]
        \item \label{prop:same_edge:1}
            $\varphi_1\approx\varphi_2$, that is, the isogenies are equivalent.
        \item \label{prop:same_edge:2}
            There are elements $k_A,k_B \in K^\times $ such that $k_AI_1=I_2$, $k_BJ_1=J_2$ and $x_2k_A=k_Bx_1$.
        \item \label{prop:same_edge:3}
            There are elements $k_A,k_B \in K^\times $ such that $k_AI_1=I_2$, $k_BJ_1=J_2$ and $k_BM_1=M_2$.
    \end{enumerate}
\end{proposition}
\begin{proof}
    The equivalence of \ref{prop:same_edge:1} and \ref{prop:same_edge:2} follows from the fact that $\cF$ is a functor.
    The fact that \ref{prop:same_edge:2} implies \ref{prop:same_edge:3} is clear.
    So, assume \ref{prop:same_edge:3}.
    We get
    \[ k_Bx_1 I_1 = k_BM_1 = M_2 = x_2I_2 = x_2 k_A I_1. \]
    Let $u = \frac{x_2k_A}{k_Bx_1}$.
    Then $u$ is a unit of the order $(I_1:I_1)$.
    If we let $k_0=k_Au^{-1}$ we then obtain $k_0I_1=I_2$, $k_BJ_1=J_2$ and $x_2k_0 = k_Bx_1$, by which \ref{prop:same_edge:2} holds.
\end{proof}

\begin{corollary}\label{cor:eq_cl_toJ}
	Let $\varphi_1\colon A_1\to B$ and $\varphi_2\colon A_2\to B$ be isogenies with the same codomain.
	Let $J = \cF(B)$ and, for $i=1,2$, let $I_i = \cF(A_i) $, $x_i = \cF(\varphi_i)$ and $M_i= x_iI_i \subseteq J$.
    The following are equivalent:
    \begin{enumerate}[(i)]
        \item \label{cor:eq_cl_toJ:1}
         $\varphi_1\approx\varphi_2$, that is, the isogenies are equivalent.
        \item \label{cor:eq_cl_toJ:2}
        There exists an element $k_A\in K^\times$ such that $k_AI_1=I_2$ and a unit $u\in (J:J)^\times$ such that $x_2k_A = x_1 u$.
        \item \label{cor:eq_cl_toJ:3}
        There exists a unit $u\in (J:J)^\times$ such that $uM_1=M_2$.
    \end{enumerate}
    Moreover, if $A_1=A_2$ and we let $\cF(A_1)=\cF(A_2)=I$ then \ref{cor:eq_cl_toJ:1}-\ref{cor:eq_cl_toJ:3} above are also equivalent to
    \begin{enumerate}[(i)]
        \setcounter{enumi}{3}
        \item \label{cor:eq_cl_toJ:4}
        there are units $v\in (I:I)^\times$ and $u\in (J:J)^\times$ such that $x_2v=ux_1$.
    \end{enumerate}
\end{corollary}
\begin{proof}
    The equivalence of \ref{cor:eq_cl_toJ:1} and \ref{cor:eq_cl_toJ:2} follows by Proposition \ref{prop:same_edge}.
    Assume \ref{cor:eq_cl_toJ:1}.
    By Proposition~\ref{prop:same_edge}, there is an element $u\in K^\times$ such that $uJ=J$ and $uM_1=M_2$.
    In particular, $u\in (J:J)^\times$, and~\ref{cor:eq_cl_toJ:3} holds.
    Assume now \ref{cor:eq_cl_toJ:3}, so that $ux_1 I_1 = x_2I_2$.
    By setting $k_A=ux_1/x_2$ and $k_B=u$, we see that Proposition~\ref{prop:same_edge}.\ref{prop:same_edge:3} holds, which is equivalent to \ref{cor:eq_cl_toJ:1}.
    Finally if $A_1=A_2$ the equivalence of \ref{cor:eq_cl_toJ:2} and \ref{cor:eq_cl_toJ:4} is clear.
\end{proof}

For ordinary elliptic curves with non-maximal endomorphism ring, prime degree isogenies are always vertical \cite[\S 2.7]{Sutherland_volcanoes}; this is another fact that fails in general:

\begin{example}\label{ex:hor_instead_of_vert}
	Consider the isogeny class with label \avlink{2.11.a\_ac} given by the Weil polynomial $h(x)=x^4-2x^2+121$, and as usual denote its Frobenius order by $R$.
	Then $R$ has a unique singular maximal ideal $\frm$ of index $[R:\frm]=2$.
    The $4$ overorders of $R$ form a chain of inclusions $R\subset S_4 \subset S_2 \subset \cO_K$ all of index $2$.

	A computation shows there are two weak equivalence classes of fractional $R$-ideals with multiplicator ring $S_4$
 	: $\cW_{S_4}(R) = \lbrace [S_4], [ S_4^\dagger] \rbrace$.
Set $J=S_4^\dagger$, which is not invertible as a fractional $S_4$-ideal.
	A further computation shows that $J$ contains exactly $3$ fractional $R$-ideals $I_1,I_2,I_3 \subset J$ with $[J:I_i]=2$.
	For each such $I_i$, the inclusion $1 \cdot I_i\subset J$ represents a minimal isogeny of degree $2$.
	One of these fractional ideals is the trace dual ideal, say $I_1=S_2^\dagger$; since $S_4 \subsetneq S_2$ this implies the isogeny  corresponding to $1 \cdot I_1 \subset J$ is a minimal descending isogeny.

	The other two ideals $I_2,I_3$ are non-isomorphic invertible fractional $S_4$-ideals.
	In particular, the isogenies induced by $1 \cdot I_2 \subset J$ and $1 \cdot I_3 \subset J$ are horizontal.
	For both of these isogenies the source and target are not weakly equivalent as ideals.
	Moreover, using Proposition~\ref{prop:same_edge}, one sees that these isogenies are not equivalent.
\end{example}

With these results, we are ready to compute isogeny equivalences.

\begin{algorithm}\label{alg:equiv_isog}$ $

    Input: For $R$ the Frobenius order of an isogeny class, inclusions $x_1I_1\subset J_1$ and $x_2I_2\subset J_2$ of $R$-ideals representing two isogenies.

    Output: Whether the two isogenies are equivalent.

    \begin{enumerate}
        \item If $I_1$ and $I_2$, or $J_1$ and $J_2$, are not isomorphic then return \emph{false}.
        \item Compute the multiplicator rings $S$ of $I_1$ and $T$ of $J_1$.
        \item Compute the unit groups $S^\times, T^\times$ as subgroups of $\cO_K^\times$.
        \item Compute elements $y,z\in K^\times$ such that $I_1=yI_2$ and $J_1=zJ_2$.
        \item If $\frac{x_1y}{x_2z}\in S^\times T^\times$ then return \emph{true}; else return \emph{false}.
    \end{enumerate}
\end{algorithm}
\begin{proof}
    The algorithm is correct by Corollary~\ref{cor:eq_cl_toJ}.\ref{cor:eq_cl_toJ:4}.\\
\end{proof}

We caution that composition is not well defined on equivalence classes:

\begin{proposition}\label{prop:multivalue}
	Let $\varphi\colon A\to B$ and $\psi\colon B\to C$ be isogenies and $\rho$ an automorphism of $B$; then $\varphi$ and $\rho\circ\varphi$ are equivalent isogenies.
	Let $H = \cF(A)$, $I = \cF(B)$, $J = \cF(C)$, $x = \cF(\varphi)$, $y = \cF(\psi)$ and $z=\cF(\rho)\in (I:I)^\times$.
    Then the compositions $\psi\circ\varphi$ and $\psi\circ\rho\circ\varphi$ are equivalent if and only if $z \in (H:H)^\times (J:J)^\times$.
\end{proposition}
\begin{proof}
    The compositions $\psi\circ\varphi$ and $\psi\circ\rho\circ\varphi$ correspond to the inclusions
    $yxH\subset J$ and $yzxH\subset J$, respectively.
    By corollary~\ref{cor:eq_cl_toJ}.\ref{cor:eq_cl_toJ:4}, the compositions are equivalent if and only if there are units $u\in (H:H)^\times$ and $v\in (J:J)^\times$ such that $uyx = vyzx$.
    This happens precisely when $z = uv^{-1}$, thus completing the proof.
\end{proof}

Finally we end this section with a result which will be useful below:

\begin{lemma}\label{lem:bijection}
	Let $R$ be the Frobenius order of an isogeny class and $I$ and $J$ be fractional $R$-ideals.
	Assume that $I$ and $J$ are weakly equivalent, so that $(I:J)$ is an invertible fractional $(I:I)$-ideal.
	Let $L$ be an invertible fractional $R$-ideal such that $L(I:I)=(I:J)$, which exists by Lemma~\ref{lem:pic_ext}.
    Then the following maps are inverses of each other:
    \begin{align*}
        \{ \text{fractional $R$-ideals $M\subseteq J$} \} & \longrightarrow \{ \text{fractional $R$-ideals $N \subseteq I$}  \} \\
        M & \longmapsto ML \\
        NL^{-1} & \mathrel{\reflectbox{\ensuremath{\longmapsto}}} N.
    \end{align*}
    This bijection preserves inclusions and indices: $[J:M]=[I:ML]$.
\end{lemma}
\begin{proof}
    Since $L$ is invertible, clearly the map is injective and inclusion preserving. 
    Let $p$ be a rational prime dividing $[J:M]$ and let $x$ be a generator of $L_p$ over $R_p$.
    Then $I_p = x_pJ_p$ and $(ML)_p = x_pM_p$ which shows that the bijection preserves indices. 
    To prove that the map is surjective, choose $N\subseteq I$, let $L^{-1}=(R:L)$, and denote by $S$ the multiplicator ring of $I$.
    Since $L^{-1}S\cdot LS=R S = S$ and the inverse is unique, we get that $L^{-1} S = (LS)^{-1}$, which implies that $L^{-1}S = (J:I)$.
    Now consider the fractional $R$-ideal $M=L^{-1}N$; we have $ML=N$.
    To show that $M$ is a preimage of $N$, we need to show that $M\subseteq J$, but this is the case since
    \[ M \subseteq MS = L^{-1}NS = L^{-1}S NS = (J:I)N \subseteq (J:I)I=J. \]
\end{proof}

\section{Isogeny graphs}\label{sec:graph}

We finally define the main object of study of this article:

\begin{definition}\label{def:graphs}\
   \begin{enumerate}
       \item The \emph{full isogeny graph} $\cG$ of an isogeny class $\cI_h$ is the directed multigraph whose vertices $\cV$ are isomorphism classes of abelian varieties in $\cI_h$ and whose edges $\cE$ are equivalence classes of isogenies of degree $>1$.        
       \item Let $D$ be a positive integer.
       The \emph{$D$-isogeny graph} $\cG^D$ is the subgraph of $\cG$ with edges $\cE^D$ of degree dividing $D$.
       \item Let $D$ and $r$ be positive integers. The graph $\cG^{D,r}$ is the subgraph of $\cG$ with edges $\cE^{D,r}$ of degree dividing $D$ and length at most $r$.
   \end{enumerate}
\end{definition}

\begin{proposition}\label{prop:missing_edges}
    Let $R$ be the Frobenius order of an isogeny class $\cI_h$.
    A minimal isogeny of degree $\ell^n$ exists in $\cI_h$ if and only if $R$ has a maximal ideal of index $\ell^n$.

    In particular, if there is no maximal ideal of $R$ of index dividing $D$ then the set $\cE^D$ of edges of $\cG^D$ is empty.
\end{proposition}
\begin{proof}
    Let $\varphi\colon A\to B$ be a minimal isogeny of degree $\ell^n$ in $\cI_h$.
    Then the existence of a maximal $R$-ideal $\frl$ with index $\ell^n$ follows from Proposition~\ref{prop:isog} parts \ref{prop:isog:deg} and \ref{prop:isog:min} since any finite simple $R$-module is isomorphic to $R/\frl$ for some maximal ideal $\frl$.

    Conversely, given a maximal $R$-ideal $\frl$ of index $\ell^n$, consider the inclusion $1\cdot \frl\subset R$.
    Let $\varphi\colon A\to B$ be an isogeny in $\cI_h$ such that $\frl = \cF(A)$, $R = \cF(B)$ and $1 = \cF(\varphi)$.
    Then $\varphi$ has degree $\ell^n$ by Proposition~\ref{prop:isog}.\ref{prop:isog:deg}.
\end{proof}

\begin{remark}\label{rmk:field_ext}
    Deligne's functor $\cF$ is particularly well suited for studying base field extensions; see \cite{MarBaseExt} for a comprehensive treatment.
    We collect here a few facts which will be used in the next examples.
    Let $\cI_h$ be an ordinary isogeny class defined over $\F_q$ with squarefree Weil polynomial $h(x)$, and $R_1 = \ZFV$ be its Frobenius order.
    For $i$ a positive integer, denote by $h_i(x)$ the Weil polynomial of the isogeny class $\cI_h \otimes \F_{q^i}$, and we assume further that $h_i(x)$ is also squarefree.
    The Frobenius order $R_i$ of $\cI_{h_i}$ can be identified with the order $\Z[\pi^i,\overline{\pi}^i]$ in $K$.
    Denote by $\cG$ and $\cG_i$ the full isogeny graphs of $\cI_h$ and $\cI_{h_i}$, respectively.

    If $A$ is in $\cI_{h_i}$ and $I$ is the fractional $R_i$-ideal corresponding to $A$, then $A$ can be defined over $\F_q$ if and only if $R_1\subseteq (I:I)$.
    More generally, the natural inclusion $\ICM(R_1) \subseteq \ICM(R_i)$ implies that $\cG$ is a subgraph of $\cG_i$.
    Furthermore, for any $A,B \in \cI_{h}$, we have $\Hom_{\F_q}(A,B)=\Hom_{\F_{q^i}}(A,B)$, from which we conclude that the sets of $\F_q$-edges and of $\F_{q^i}$-edges between two vertices of $\cG$ coincide.
    Analogous statements hold for the $D$-isogeny graph $\cG^D$.

    Finally, if $\cI_h$ is geometrically simple and ordinary then $h_i(x)$ is an irreducible polynomial for every $i>1$.
    Therefore, the considerations above apply to the extension to~$\overline{\F}_q$.
\end{remark}

\begin{corollary}\label{cor:index-divisibility}
  Let $\ell$ be a prime number and $\cG$ be the $\ell$-isogeny graph of an isogeny class defined over $\F_q$ with Frobenius order $R_1$.
	Let $\cG'$ be the $\ell$-isogeny graph for the same isogeny class base changed $\F_{q^i}$ for $i>1$, and with Frobenius order $R_i$.
	In this case if $\cG$ has no edges and $\cG'$ has edges then $\ell$ divides $[R_1: R_i]$.
\end{corollary}
\begin{proof}
  Proposition~\ref{prop:missing_edges} implies that $R_i$ has a maximal ideal $\frm$ of index $\ell$ while $R_1$ does not.
Note that $R_1/\frm R_1$ is an $R_i / \frm$-vector space of dimension $n > 1$, since $R_1$ does not have a maximal ideal of index $\ell$.
Consider the map of finite abelian groups $\sigma\colon R_1/R_i \to (R_1/\frm R_1)/(R_i/\frm)$.
We have that $[R_1: R_i] = \ell^{n-1} \cdot \vert \ker(\sigma) \vert$ which implies that $\ell$ divides $[R_1: R_i]$.
\end{proof}

To illustrate Remark~\ref{rmk:field_ext}, we present an example of a $2$-isogeny graph of abelian surfaces over $\F_3$ with no edges.
We explain how $2$-isogenies appear after base change to $\F_{3^i}$ for appropriate $i>1$.
Note that there exist arbitrarily large $i$ where the $2$-isogeny graph over $\F_{3^i}$ has no edge.

\begin{example}\label{ex:no-edges-gain_edges}
    Consider the isogeny class $\cI_h$ of geometrically simple abelian surfaces over $\F_3$ with LMFDB label \avlink{2.3.ad\_f} and Weil polynomial
    $h(x)=x^4 - 3 x^3 + 5 x^2 - 9 x + 9$.
    Over the fields $\F_3,\F_9,\F_{81},$ and $ \F_{6561}$ the $2$-isogeny graph has 1, 5, 145 and  672945 vertices respectively, all with no edges.

   The Frobenius order $R_1$ of $\cI_h$ is the maximal order $\cO_K$; it has a unique maximal ideal $\frm$ above $2$ which is of index $16$.
   By Proposition~\ref{prop:missing_edges}, the $2$-isogeny graph of $\cI_h$ therefore has no edges. 

   Now for $i$ a positive integer and $\cI_{h_i}$ the base change of $\cI_h$ to $\F_{3^i}$, let $R_i$ be the Frobenius order of $\cI_{h_i}$.
   Corollary~\ref{cor:index-divisibility} implies that if $[R_1:R_i]$ is coprime to $2$, there will be no edges in the $2$-isogeny graph over $\F_{3^i}$.
 Table~\ref{tab:table-of-indices} records first the cardinality of $\ICM(R_i)$ -- which we observe grows very quickly -- 
 and the indices of the maximal primes above $2$.
  	Then, to track the development of $2$-isogenies we define $[R_1:R_i]_\mathrm{new}$ to be the index $[R_1:R_i]$ divided by the least common multiple of $[R_1:R_j]$ for proper divisors $j$ of $i$.

\begin{table}[!htbp]
\begin{center}
	\begin{tabular}{|c|c|c|l|l|}
        \hline $i$ & $\#\ICM(R_i)$ &$\left[ [R_i:\frm] : 2 \in \frm \right] $ & $[R_1:R_i]_\mathrm{new}$ & $[R_1:R_i] / [R_1:R_i]_{\mathrm{new}}$ \\ \hline \hline
		$1$ & $1$ & $[ 16 ]$ & $ 1$ & $1$ \\
		$2$ & $5$ & $[ 16 ]$ & $ 3^{2}$ & $1$ \\
		$3$ & $25$ & $[ 16 ]$ & $ 23$ & $1$ \\
		$4$ & $145$ & $[ 16 ]$ & $ 29$ & $3^{2}$ \\
		$5$ & $5246$ & $[ 2 ]$ & $ 2^{6}\cdot 61$ & $1$ \\
		$6$ & $8075$ & $[ 16 ]$ & $ 3^{2}\cdot 17$ & $3^2\cdot 23$ \\
		$7$ & $201707$ & $[ 16 ]$ & $ 29^{2}\cdot 1667$ & $1$ \\
		$8$ & $672945$ & $[ 16 ]$ & $ 23^{2}\cdot 103$ & $3^{2}\cdot 29$ \\
		$9$ & $234025$ & $[ 16 ]$ & $ 17^{2}\cdot 251$ & $23$ \\
		$10$ & $10677440$ & $[ 2 ]$ & $ 2^{8}$ & $2^{6}\cdot 3^{2}\cdot 61$ \\
		$11$ & $10716555$ & $[ 16 ]$ & $ 419^{2}\cdot 25453$ & $1$ \\
		$12$ & $3464756875$ & $[ 16 ]$ & $ 23\cdot 61^{2}$ & $3^{4}\cdot 17\cdot 23\cdot 29$ \\
		$13$ & $-$ & $[ 16 ]$ & $ 5^{2}\cdot 13^{6}\cdot 131$ & $1$ \\
		$14$ & $-$ & $[ 16 ]$ & $ 2731$ & $3^{2}\cdot 29^{2}\cdot 1667$ \\
		$15$ & $-$ & $[ 2 ]$ & $ 569^{2}\cdot 719$ & $2^{6}\cdot 23\cdot 61$ \\
		$16$ & $-$ & $[ 16 ]$ & $ 15473$ & $3^{2}\cdot 23^{2}\cdot 29\cdot 103$ \\
		$17$ & $-$ & $[ 16 ]$ & $ 9283^{2}\cdot 10449287$ & $1$ \\
		$18$ & $-$ & $[ 16 ]$ & $ 3^{2}\cdot 17^{4}\cdot 179$ & $3^{4}\cdot 17^{2}\cdot 23\cdot 251$ \\
		$19$ & $-$ & $[ 16 ]$ & $ 191\cdot 1901^{2}\cdot 276337$ & $1$ \\
		$20$ & $-$ & $[ 2 ]$ & $ 2^{6}\cdot 199^{2}\cdot 1301$ & $2^{14}\cdot 3^{2}\cdot 29\cdot 61$ \\\hline
	\end{tabular}
\end{center}
\caption{Information about the isogeny class \avlink{2.3.ad\_f}. See Example~\ref{ex:no-edges-gain_edges} for a description.}
\label{tab:table-of-indices}
\end{table}

\end{example}

\begin{theorem}\label{thm:conn}
Let $D>1$, $R$ be the Frobenius order of the isogeny class $\cI_h$, and $\cS^\mathrm{sing}$ be the set of maximal $R$-ideals containing the conductor $\frf_R=(R:\cO_K)$ of $R$.\

\begin{enumerate}[(i)]
\item \label{thm:conn:conn}
    Suppose that there exists a set $\cS$ of maximal ideals of $R$ containing a generating set of $\Pic(R)$ such that the integer
\[ N_{\cS}=\mathrm{lcm}\left(
\left\{ \vert R/\frm\vert : \frm\in \cS^\mathrm{sing} \cup \cS \right\}
\right). \]
divides $D$.  Then the isogeny graph $\cG^D$ is connected (as a directed graph).
\item \label{thm:conn:diam}
    Let $\cS$ be a set of maximal $R$-ideals satisfying the conditions in \ref{thm:conn:conn}.
    Then the diameter of $\cG^{D,1}$ satisfies:
$$ \mathrm{diam}(\cG^{D,1})\leq
\begin{cases}
	2\left(\len_R(\cO_K/R)-1\right) + \sum_{\mathfrak{l}\in \cS}\left(\mathrm{ord}_{\Pic(R)}([\mathfrak{l}])-1\right) & \mbox{ if $R$ is Bass,} \\
	2\left(\len_R(\cO_K/\mathfrak{f}_R)-1\right) + \sum_{\mathfrak{l}\in \cS}\left(\mathrm{ord}_{\Pic(R)}([\mathfrak{l}])-1\right) & \mbox{ otherwise,}
\end{cases}$$
where $\mathrm{ord}_{G}(h)$ for an element $h$ of a group $G$ denotes the order of $h$ in $G$.
\end{enumerate}
\end{theorem}
\begin{proof}

Let $\cS=\{ \mathfrak{l}_1,\ldots,\mathfrak{l}_n \}$ and $N_{\cS}$ be as in the statement.
We first show that any two vertices in the same weak equivalence class are connected by a horizontal edge in $\cE^D$, then that the weak equivalence classes are connected.
This construction also gives the bounds of \ref{thm:conn:diam}.

Let $I$ and $J$ be fractional $R$-ideals in the same weak equivalence class.
By Proposition~\ref{prop:pic_weak}, this means that $I=LJ$ for some invertible fractional $R$-ideal $L$.
Since $\cS$ generates $\Pic(R)$, there exists an element $a\in K^\times$  such that
\[ aI = \mathfrak{l}_1^{e_1}\cdots\mathfrak{l}_n^{e_n}J, \]
for integers $0\leq e_i < \mathrm{ord}_{\Pic(R)}([\mathfrak{l}_i])$.
Hence, $\len_R\left(J/aI\right) = \sum_{i=1}^n e_i$,
and $a$ represents an edge which is a composition of $e_i$ minimal edges of degree $\vert R/\mathfrak{l}_i \vert$ by Proposition~\ref{prop:isog}, which each belong to $\cE^D$ since $N_\cS$ divides $D$ by assumption.

Now, consider any two fractional $R$-ideals $I$ and $J$.
By \cite[Remark~5.4]{MarICM18}, there are invertible $R$-ideals $L_I$ and $L_J$ such that
$\mathfrak{f}_R L_I \subseteq I \subseteq \cO_K L_I$ and $\mathfrak{f}_R L_J \subseteq J \subseteq \cO_K L_J$.
The ideals $\cO_K L_I$ and $\mathfrak{f}_R L_J$ are weakly equivalent, since they both have multiplicator ring $\cO_K$ and there is only one such weak equivalence class of ideals as every $\cO_K$-ideal is invertible. 
Therefore an edge $[\cO_K L_I]\to [\frf_R L_J]$ can be obtained by composing at most $\sum_{i=1}^n \left(\mathrm{ord}_{\Pic(R)}([\frl_i])-1\right)$ minimal edges.

We have isomorphisms of finite $R$-modules $\cO_K/\mathfrak{f}_R \simeq \cO_K L_I/\mathfrak{f}_R L_I \simeq \cO_K L_J/\mathfrak{f}_R L_J$.
The simple factors of $\cO_K/\mathfrak{f}_R$ are $R/\mathfrak{m}$ for $\mathfrak{m}\in \cS^\mathrm{sing}$.
By Proposition~\ref{prop:isog}, each such simple factor is equivalent to an edge of degree $\vert R/\mathfrak{m}\vert$ which belongs to $\cE^D$ by the assumption on $N_\cS$.
This concludes the proof that $\cG^D$ is connected.

To obtain \ref{thm:conn:diam}, we now note that the paths $[I]\to [\cO_K L_I]$ and $[\mathfrak{f}_R L_J]\to [J]$ in $\cG^{D,1}$ have length bounded by $\len_R\left(\cO_K/\mathfrak{f}_R\right) -1$, where the $-1$ comes from the fact that $\mathfrak{f}_R$ and $\cO_K$ are in the same weak equivalence class.
When $R$ is Bass, the weak equivalence classes are represented by the overorders of $R$, so we may consider the composition series of $\cO_K/R$ instead of $\cO_K/\mathfrak{f}_R$.  This yields an improved upper bound and concludes the proof of \ref{thm:conn:diam}.
\end{proof}

\section{Computing the isogeny graph}\label{sec:algorithms_GD}

We now give an algorithm to compute the $D$-isogeny graph $\cG^D$ of the isogeny class $\cI_h$ with Frobenius order $R$.
Throughout, for each weak equivalence class $w$ of fractional $R$-ideals, we fix a representative $W_w$ and write $\frw_w$ for its class in $\ICM(R)$.
Further, we represent each ideal class $[M]\in\ICM(R)$ by a unique pair $(w_M,\fra_M)$ where $w_M$ is the weak equivalence class of $M$ and $\fra_M=[(M:W_{w_M})] \in \Pic((M:M))$.
Finally, for each overorder $S$ of $R$ let $e_S$ be the map of Lemma \ref{lem:pic_ext}, and for each class in $\Pic(S)$, fix a preimage $\fra\in \Pic(R)$ via $e_S$ and an invertible fractional $R$-ideal $I_\fra$ representing $\fra$.
We denote the set of chosen preimages in $\Pic(R)$ by $\cP_S$.

The computation of $\cG^D$ is divided in two steps:
First we compute the subgraph $\cG^{D,1}$ containing only minimal edges in Algorithm~\ref{alg:min_isog}, and then we recover $\cG^D$ by taking all possible compositions in Algorithm~\ref{alg:compositions}.

\begin{algorithm}\label{alg:min_isog}$ $

    Input: An isogeny class of ordinary abelian varieties over a finite field $\F_q$ given by its Frobenius order $R$, and an integer $D>1$.

    Output: The subgraph $\cG^{D,1}$ of the $D$-isogeny graph $\cG^D$.
    \begin{enumerate}
        \item Compute the set $\cV$ of vertices $(t,e_T(\frb))$ by computing all weak equivalence classes $t$ of $R$ and ideal classes $\frb\in\cP_T$, for $T$ the multiplicator ring of $t$.
        \item For each weak equivalence class $t$, let $T$ be the multiplicator ring of $t$ and:
        \begin{enumerate}
            \item\label{alg:min_isog:3b} Compute the fractional $R$-ideals $M$ contained in $W_t$ such that the inclusion $M\subset W_t$ is maximal and the index $[W_t:M]$ divides $D$.
            \item\label{alg:min_isog:3c} Compute the orbits of these submodules under the action of $T^\times$, and choose a representative $M$ for each orbit.
            \item\label{alg:min_isog:ED1} For each orbit representative $M$ do:
            \begin{enumerate}
                \item \label{alg:min_isog:2di}
                Compute the vertex $(s,e_S(\fra))\in \cV$ that represents the ideal class of $M$, where $S=(M:M)$ and $\fra\in\cP_S$, together with $x\in K^\times$ such that $x W_sI_\fra = M$.
                \item\label{alg:min_isog:min_edges} For every $\frb\in\cP_T$ do:
                \begin{enumerate}
                    \item Compute $y\in K^\times$ such that $yI_{\fra'}S=I_\fra I_\frb S$ where $\fra'$ is the unique element in $\cP_S$ satisfying $e_S(\mathfrak{ab}) = e_S(\fra')$.
                    \item Yield the minimal edge $(s,e_S(\fra'))\to (t,e_T(\frb))$ with label $xyW_sI_{\fra'}\subset W_tI_\frb$.
                \end{enumerate}
            \end{enumerate}
        \end{enumerate}
    \end{enumerate}
\end{algorithm}
\begin{proof}
    The set of vertices $\cV$ is complete by Proposition~\ref{prop:pic_weak}.

    Let $I$ and $J$ be fractional $R$-ideals and $a \in K^\times$ be such that $aI \subset J$ is maximal with index dividing $D$.  To show that the set $\cE^{D,1}$ of minimal edges returned by the algorithm is complete, we prove that there is an edge in $\cE^{D,1}$ representing the equivalence class of the corresponding minimal isogeny.
    Let $t$ be the weak equivalence class of $J$ and $\frb\in \cP_T$ be such that $(J:W_t)\in e_T(\frb)$.
    Thus there exists $b\in K^\times$ such that $bI_\frb T = (J:W_t)$, where $I_\frb$ is the chosen representative of $\frb$.
    Consider the invertible fractional $R$-ideal $L=bI_\frb$.
    By construction, we have $LW_t = J$.
    Set $M=L^{-1}(aI)$.
    By Lemma~\ref{lem:bijection}, $M\subset W_t$ is a maximal inclusion of index dividing $D$, that is, $M$ is one of the fractional $R$-ideals computed in step \ref{alg:min_isog:3b}.
    If $u$ is a unit of $T = (W_t:W_t)=(J:J)$ then the isogenies $aI\subset J$ and $uaI\subset J$ are equivalent by Corollary~\ref{cor:eq_cl_toJ}.\ref{cor:eq_cl_toJ:4}.
    Hence, we can assume that $M$ is the orbit representative picked in step \ref{alg:min_isog:ED1}.
    We have $MI_\frb = (1/b)aI$ and $W_tI_\frb=(1/b)J$.
    By Proposition~\ref{prop:same_edge}.\ref{prop:same_edge:1}, the isogeny $aI\subset J$ is equivalent to $MI_\frb\subset W_tI_\frb$.
    Finally, the isogeny $M I_\frb \subset W_t I_\frb$ is equivalent to $xyW_sI_{\fra'} \subset W_t I_\frb$ constructed in step \ref{alg:min_isog:min_edges} by Corollary~\ref{cor:eq_cl_toJ}.\ref{cor:eq_cl_toJ:3}.
    Hence every minimal isogeny with index dividing $D$ is represented in $\cE^{D,1}$,
    and $\cE^{D,1}$ represents all minimal isogenies of degree dividing $D$ by Proposition~\ref{prop:isog}.\ref{prop:isog:incl}.

    We now prove that distinct elements of $\cE^{D,1}$ represent inequivalent isogenies.
    If the source or target vertices of two edges are not the same, the edges represent inequivalent isogenies.
    Therefore we consider two edges produced in step \ref{alg:min_isog:min_edges} with the same source and with target $(t,e_T(\frb))$, where $T$ is the multiplicator ring of $t$ and $\frb\in \cP_T$.
    A priori, we know that the sources of the two edges are of the form $(s_i,e_S(\fra_i\frb))$ for $i=1,2$.
    Since we are assuming that they are equal, we get $s_1=s_2$ and $\fra_1=\fra_2$, which we then denote simply by $s$ and $\fra$, respectively.
    Because the edges are distinct, they come from fractional $R$-ideals $M_1,$ and $M_2$ contained in $W_t$ which do not belong to the same orbit under $T^\times$.
    Write $M_1=x_1W_sI_\fra$ and $M_2=x_2W_sI_\fra$ for elements $x_1,x_2\in K^\times$ so that the edges represent the inclusions $M_1I_\frb\subseteq W_tI_\frb$ and $M_2I_\frb\subseteq W_tI_\frb$.
    Corollary~\ref{cor:eq_cl_toJ}.\ref{cor:eq_cl_toJ:3} then implies that these edges represent inequivalent isogenies.
\end{proof}

\begin{algorithm}\label{alg:compositions}$ $

    Input: An isogeny class of ordinary abelian varieties over a finite field $\F_q$ given by its Frobenius order $R$, and an integer $D>1$.

    Output: The isogeny graph $\cG^D$.
    \begin{enumerate}
        \item Run Algorithm~\ref{alg:min_isog} to compute the vertices $\cV$ and minimal edges $\cE^{D,1}$ of $\cG^D$.   Write $\cE^{=n}$ for the edges of degree exactly $n$, $\cE^{=n,1}$ for the minimal edges of degree $n$, and initialize $\cE^{=n} = \cE^{=n,1}$ for $1 < n \mid D$.
        \item Loop over the divisors $n\neq 1$ of $D$ in increasing order, and
            for each divisor $d_2>1$ of $n$, loop over pairs $(E_2,E_1) \in \cE^{=d_2,1}\times \cE^{=n/d_2}$ such that
            the source of the minimal edge $E_2\colon (w_2,\fra_2)\to (w_3,\fra_3)$ equals the target of the edge $E_1\colon (w_1,\fra_1)\to (w_2,\fra_2)$.
            \begin{enumerate}
                \item Compute a transversal $\cU$ of the quotient $\cO_2^\times/(\cO_1^\times\cO_3^\times \cap \cO_2^\times)$, where $\cO_i$ is the multiplicator ring of $w_i$.
                \item\label{alg:compositions:U} For each $u\in \cU$ do:
                \begin{enumerate}
                        \item Compute the degree-$n$ composition $E\colon (w_1,\fra_1) \to (w_2,\fra_2) \to (w_3,\fra_3)$ labelled by $x_1ux_2W_{w_1}I_{\fra_1} \subset W_{w_3}I_{\fra_3}$.
                        \item If Algorithm~\ref{alg:equiv_isog} returns that there is no already computed edge equivalent to $E$ in $\cE^{=n}$, add $E$ to $\cE^{=n}$.
                \end{enumerate}
            \end{enumerate}
        \item Construct $\cG^D$ from $\cV$ and the edges $\cE^{=n}$ for $1 < n \mid D$ and return it.
    \end{enumerate}
\end{algorithm}
\begin{proof}
    The labels of the minimal edges in $\cE^{D,1}$, and consequently of all of $\cE^D$, are chosen so that the ideals representing source and target are equal whenever source and target coincide.
    In other words, we do not need to consider isomorphism between different representatives when computing the compositions.
    Let $n>1$ be a divisor of $D$.
    Every edge of degree $n$ can be written as a minimal edge $E_2$ of degree $d_2|n$ following, if $d_2<n$, an edge $E_1$ of degree $n/d_2<n$.
    Since we are looping over the divisors $n$ of $D$ in increasing order, all such edges $E_1$ have already been computed.
    In step \ref{alg:compositions:U}, we loop over the transversal $\cU$ to consider all the contributions coming from the automorphisms of the source of $E_2$ given by Proposition~\ref{prop:multivalue},
    and therefore the set $\cE^{=n}$ contains representatives of all equivalence classes of the edges of degree $n$.
    The use of Algorithm~\ref{alg:equiv_isog} and the fact that the edges produced for each $u \in \cU$ are all pairwise not equivalent by Corollary~\ref{cor:eq_cl_toJ}.\ref{cor:eq_cl_toJ:4} guarantee that $\cE^{=n}$ contains at most one representative for each equivalence class.
\end{proof}

\section{The Picard group action}\label{sec:alg_storing}

Algorithms \ref{alg:min_isog} and \ref{alg:compositions}  will struggle when $\Pic(R) = \Pic(\ZFV)$ is large.
In this situation, we can try to take advantage of additional structure to both improve the runtime and reduce the size of the output.
Indeed, the group $\Pic(R)$ acts on both the vertices and edges of $\cG$ via ideal multiplication.
While this action is not free in general, we can use it to define free actions of specific quotient groups on subsets of the edges in $\cG$.

\begin{definition} Suppose that $s$ and $t$ are weak equivalence classes of $R$ with multiplicator rings $S$ and $T$, $e_S \colon  \Pic(R) \to \Pic(S)$ and $e_T \colon  \Pic(R) \to \Pic(T)$.
\begin{enumerate}
\item Write $\dom \colon  \cE \to \cV$ for the map that sends each isogeny to its domain, and $\cod \colon  \cE \to \cV$ for the map that sends each isogeny to its codomain.
\item Set $\cE_{s, t} = \{\varphi \in \cE : \dom(\varphi) \in s, \cod(\varphi) \in t\}$, and write $\cE_{s,t}^D$ for the set of elements of $\cE_{s,t}$ of degree dividing $D$ and $\cE_{s,t}^{D,r}$ for the set of elements of $\cE_{s,t}^D$ of length $r$.
\item Finally, define the quotient group $G_{S, T} = \Pic(R) / (\ker(e_S) \cap \ker(e_T))$.
\end{enumerate}
\end{definition}

We now extend \cite[Theorem 4.6]{MarICM18} from vertices to edges:

\begin{proposition}
Let $A$ and $B$ belong to an isogeny class with isogeny graph $\cG$ and Frobenius order $R$, $I = \cF(A)$, $J = \cF(B)$, $\varphi \colon  A \to B$ be an isogeny, $a = \cF(\varphi)$, 
and $M$ be an invertible fractional $R$-ideal with class $[M] \in \Pic(R)$.
\begin{enumerate}
\item Ideal multiplication induces an action of $\Pic(R)$ on $\cG$: we denote by $[M] \cdot \varphi$ the isogeny represented by the inclusion $aMI \subseteq MJ$.
\item The action of $\Pic(R)$ on $\cG$ preserves degree and length, and thus induces actions on $\cG^D$ and $\cG^{D,r}$.
\item This action induces an action of $G_{S,T}$ on $\cE_{s,t}$, $\cE_{s,t}^D$ and $\cE_{s,t}^{D,r}$.
\end{enumerate}
\end{proposition}
\begin{proof}
The first two parts follow from Proposition \ref{prop:isog} and Lemma \ref{lem:bijection}.

To see that the action on $\cE_{s,t}$ factors through $G_{S,T}$, we show that $\ker(e_S) \cap \ker(e_T)$ acts trivially.  Let $[H] \in \ker(e_S) \cap \ker(e_T)$, $I \in s$, $J \in t$, and $xI \subseteq J$ represent an edge in $\cE_{s,t}$.  Acting by $[H]$ yields $xHI \subseteq HJ$.  Since $HS = y_1S$ and $HT = y_2T$ for $y_1,y_2 \in K$, the two inclusions represent equivalent isogenies by Proposition \ref{prop:same_edge}.
\end{proof}

\begin{lemma} \label{lem:GST_free} Let $R$ be the Frobenius order of an isogeny class, and $s$ and $t$ be weak equivalence classes of fractional $R$-ideals with multiplicator rings $S$ and $T$.
\begin{enumerate}
\item \label{lem:GST_free:1} The group $G_{S,T}$ acts freely on $\cE_{s,t}$.
\item \label{lem:GST_free:2} For each orbit $O$ of this action, we have $\dom(O) = s$ and $\cod(O) = t$.
\item \label{lem:GST_free:3} If $O$ is an orbit and $\frj \in t$, the set $\{\varphi \in O : \cod(\varphi) = \frj\}$ is a torsor for $\ker(e_T) / (\ker(e_S) \cap \ker(e_T))$.  Similarly, if $\fri \in s$, the set $\{\varphi \in O : \dom(\varphi) = \fri\}$ is a torsor for $\ker(e_S) / (\ker(e_S) \cap \ker(e_T))$.
\end{enumerate}
\end{lemma}
\begin{proof}
We may assume that $\cE_{s,t}$ is nonempty since otherwise the action is vacuously free and there are no orbits.
Suppose that $\fra \in \Pic(R)$ fixes some $\varphi \in \cE_{s,t}$.  Then it certainly fixes $\dom(\varphi)$, and thus $e_S(\fra) = 1$ by Proposition~\ref{prop:pic_weak}. Similarly, it fixes $\cod(\varphi)$ so $e_T(\fra) = 1$.  Thus the action of $G_{S,T}$ is free.

Since $e_S$ is surjective, $\Pic(R)$ acts transitively on $s$ and thus $\dom(O) = s$ for any orbit $O$ of the action of $G_{S,T}$ on $\cE_{s,t}$. Similarly, $e_T$ is surjective so $\cod(O) = t$.

Since $G_{S,T}$ acts on codomains through its projection onto $\Pic(T)$, the subgroup $\ker(e_T) / (\ker(e_S) \cap \ker(e_T))$ is precisely the set of group elements fixing the codomain. 
It thus acts on the set $\{\varphi \in O : \cod(\varphi) = \frj\}$, which is a torsor for the subgroup since $G_{S,T}$ acts freely.  The statement about $\{\varphi \in O : \dom(\varphi) = \fri\}$ follows similarly.
\end{proof}

In light of  Lemma \ref{lem:GST_free} , we modify Algorithms \ref{alg:min_isog} and \ref{alg:compositions} to return sets of representatives $\cR_{s,t}$ for the orbits of $G_{S,T}$ on $\cE_{s,t}$. This has the benefit of both dramatically shrinking the output size and reducing the amount of compositions computed when $\Pic(R)$ is large.  However, composing orbit representatives requires additional work.  We consider the subsets $\cR_{s,t}^{=n}$ and $\cR_{s,t}^{=n,1}$ consisting of edges with degree $n$ and minimal edges of degree $n$.

In this section, $s, t,$ and $u$ are weak equivalence classes of fractional $R$-ideals, with corresponding multiplicator rings $S$, $T$, and $U$.  
Furthermore, for isogenies $\varphi$ and $\varphi'$, we write $\varphi' \sim \varphi$ if $\varphi' \approx \fra \cdot \varphi$ for some $\fra \in \Pic(R)$. 
Finally, we define a multivalued composition map $\ucirc \colon  \cR_{t,u} \times \cR_{s,t} \to \cR_{s,u}$ by setting
\[\psi \ucirc \varphi = \{ \theta \in \cR_{s,u} : \exists \varphi', \psi', \theta' \text{ with } \varphi' \sim \varphi, \psi' \sim \psi, \theta' \sim \theta, \theta' = \psi' \circ \varphi'\}.\]

\begin{proposition} \label{prop:GSTcompose}
Suppose that $\fra_0 \in \Pic(R)$ satisfies $\cod(\fra_0 \cdot \varphi) = \dom(\psi)$.
Let $X$ be a set of coset representatives for $\ker(e_T) / (\ker(e_T) \cap (\ker(e_S) \cdot \ker(e_U)))$ and $Y$ be a set of coset representatives for $T^\times / (T^\times \cap (S^\times U^\times))$. Then
\begin{enumerate}
\item \label{prop:GSTcompose:part1} The $G_{S,U}$-orbits represented by $\psi \ucirc \varphi$ are the same as those represented by
\[
\psi \usq \varphi = \{\psi \circ \rho \circ (\fra_1 \fra_0 \cdot \varphi) : \fra_1 \in X, \cF(\rho) \in Y\} \subseteq \psi \ucirc \varphi.
\]
\item If $\chi = \psi \circ \rho \circ (\fra_1 \fra_0 \cdot \varphi)$ and $\chi' = \psi \circ \rho' \circ (\fra_1' \fra_0 \cdot \varphi)$ are two elements of $\psi \usq \varphi$ then $\chi \sim \chi'$ if and only if $\chi \approx \chi'$.
\end{enumerate}
\end{proposition}

\begin{proof}
By the definition of $\sim$, the orbits represented by $\psi \usq \varphi$ are contained within those represented by $\psi \ucirc \varphi$, so it suffices to check the opposite containment.  Suppose that $\varphi' \sim \varphi$ and $\psi' \sim \psi$ can be composed as $\theta' = \psi' \circ \varphi'$.  Then $\varphi' = \beta \circ (\frb_1 \cdot \varphi) \circ \alpha$ with $\cF(\alpha) \in S^\times, \cF(\beta) \in T^\times, \frb_1 \in \Pic(R)$ and $\psi' = \delta \circ (\frb_2 \cdot \psi) \circ \gamma$ with $\cF(\gamma) \in T^\times, \cF(\delta) \in U^\times, \frb_2 \in \Pic(R)$.  Since we are working with $\theta'$ up to $G_{S,U}$-equivalence, we may omit $\alpha$ and $\delta$, set $\rho = \gamma \circ \beta$, and check that each composition $(\frb_2 \cdot \psi) \circ \rho \circ (\frb_1 \cdot \varphi)$ is in the same orbit as one of the form $\psi \circ \rho \circ (\fra_1 \fra_0 \cdot \varphi)$.

Let $\fra_2 \in \ker(e_T) \cap \ker(e_U)$.  Then $\frb_2 \fra_2 \cdot (\psi \circ \rho \circ (\fra_1 \fra_0 \cdot \varphi)) \approx (\frb_2 \fra_2 \cdot \psi) \circ \rho \circ (\frb_2 \fra_2 \fra_1 \fra_0 \cdot \varphi) = (\frb_2 \cdot \psi) \circ \rho \circ (\frb_2 \fra_2 \fra_1 \fra_0 \cdot \varphi)$, so it suffices to check that $\frb_2^{-1}\frb_1 \cdot \varphi \approx \fra_2 \fra_1 \fra_0 \cdot \varphi$ for some $\fra_2 \in \ker(e_T) \cap \ker(e_U)$ and $\fra_1 \in X$.  Since the composition $(\frb_2 \cdot \psi) \circ \rho \circ (\frb_1 \cdot \varphi)$ is well defined, we must have $\cod(\frb_1 \cdot \varphi) = \dom(\frb_2 \cdot \psi)$ and thus $\cod(\frb_2^{-1} \frb_1 \cdot \varphi) = \dom(\psi) = \cod(\fra_0 \cdot \varphi)$.  But as $\fra_2$ and $\fra_1$ vary, $\fra_2\fra_1\fra_0 \cdot \varphi$ ranges over all equivalence classes of isogenies with the same codomain as $\fra_0 \cdot \varphi$. This completes the proof of part \ref{prop:GSTcompose:part1}

Now suppose that $\chi \sim \chi'$.  
To prove that $\chi \approx \chi'$, it suffices to show that $\dom(\chi) = \dom(\chi')$ and $\cod(\chi) = \cod(\chi')$ since any element of $\Pic(R)$ that acts nontrivially must change either the domain or codomain.
All elements of $\psi \usq \varphi$ have the same codomain, so suppose that $\fra_1, \fra_1' \in X$.  Since $\chi \sim \chi'$, by Lemma \ref{lem:GST_free} the domains of $\fra_1 \fra_0 \cdot \varphi$ and $\fra_1' \fra_0 \cdot \varphi$ differ by an element of $\ker(e_U)$.  By the definition of $X$, both lie in $\ker(e_T)$ as well, and thus their difference is in $\ker(e_U) \cap \ker(e_T) \subseteq (\ker(e_T) \cap (\ker(e_S) \cdot \ker(e_U)))$.  Since $X$ is a set of coset representatives, we have $\fra_1 = \fra_1'$, and thus $\dom(\chi) = \dom(\chi')$.
\end{proof}

Next, we augment Algorithm \ref{alg:equiv_isog} to give a test for when two isogenies are in the same orbit, and prove that our algorithm is correct.

\begin{algorithm}\label{alg:equiv_orbit}$ $

    Input: For $R$ the Frobenius order of an isogeny class, inclusions $x_1I_1\subset J_1$ and $x_2I_2\subset J_2$ of $R$-ideals representing two isogenies $\varphi$ and $\psi$ with weakly equivalent domains and weakly equivalent codomains.

    Output: Whether the two isogenies are in the same $\Pic(R)$-orbit.

    \begin{enumerate}
        \item Represent the ideal classes $[I_1] = (s, \fra_1'), [I_2] = (s, \fra_2'), [J_1] = (t, \frb_1'), [J_2] = (t, \frb_2')$, let $S$ be the multiplicator ring of the domains and $T$ the multiplicator ring of the codomains and choose lifts $\fra_i, \frb_i \in \Pic(R)$ of $\fra_i'$ and $\frb_i'$.
        \item \label{alg:equiv_orbit:2} Express $\fra_2 \frb_1 \fra_1^{-1} \frb_2^{-1} = \frk_S \frk_T$ for some $\frk_S \in \ker(e_S)$ and $\frk_T \in \ker(e_T)$.  If this is not possible, return \textit{false}.
        \item Compare $(\frb_2 \frb_1^{-1} \frk_T) \cdot \varphi$ to $\psi$ using Algorithm  \ref{alg:equiv_isog}.
    \end{enumerate}
\end{algorithm}
\begin{proof}
    Step \ref{alg:equiv_orbit:2} can be done by finding the preimage of $\fra_2 \frb_1 \fra_1^{-1} \frb_2^{-1}$ under the natural homomorphism $\ker(e_S) \oplus \ker(e_T) \to \Pic(R)$, and yields an element $\frb_2 \frb_1^{-1} \frk_T = \fra_2 \fra_1^{-1} \frk_S^{-1}$ whose image under $e_S$ is $\fra'_2 {\fra'}_1^{-1}$ (which maps the domain of $\varphi$ to the domain of $\psi$) and whose image under $e_T$ is $\frb'_2 {\frb'}_1^{-1}$ (which maps the codomain of $\varphi$ to the codomain of $\psi$).  If no such $\frk_S$ and $\frk_T$ exist, $\varphi$ and $\psi$ cannot be in the same orbit since it is not possible to get their domains and codomains to match.  If they do exist, then $\varphi \sim \psi$ if and only if $(\frb_2 \frb_1^{-1} \frk_T) \cdot \varphi \approx \psi$.
\end{proof}

We now adapt Algorithm \ref{alg:min_isog} to compute orbit representatives for the set of minimal isogenies.  In contrast to Algorithm \ref{alg:min_isog}, we do not need to loop over $\Pic(T)$.

\begin{algorithm}\label{alg:GST_min}$ $

    Input: An isogeny class of ordinary abelian varieties over a finite field $\F_q$ given by its Frobenius order $R$, and an integer $D>1$.

    Output: For weak equivalence classes $s,t$ and $1 < n \mid D$, the finite set $\cR_{s,t}^{=n,1}$.
    \begin{enumerate}
        \item For each weak equivalence class $t$ and $T$ the multiplicator ring of $t$, do
        \begin{enumerate}
            \item\label{alg:GST_min:3b} Compute the fractional $R$-ideals $M$ contained in $W_t$ such that the inclusion $M\subset W_t$ is maximal and the index $[W_t:M]$ divides $D$.
            \item\label{alg:GST_min:3c} Among the $M$ in the same weak equivalence class $s$, choose one representative for each $\ker(e_T) / (\ker(e_S) \cap \ker(e_T))$-orbit.  Determine the action of $T^\times$ on this quotient and pick one representative from each of these $T^\times$ orbits.
            \item\label{alg:GST_min:ED1} For each representative $M$, compute the vertex $(s,e_S(\fra))\in \cV$ that represents it with $\fra\in\cP_S$, together with $x\in K^\times$ such that $x W_sI_\mathfrak{a} = M$.  Yield the minimal edge $(s,e_S(\mathfrak{a}))\to (t,1)$ with label $xW_sI_{\mathfrak{a}}\subset W_t$.
        \end{enumerate}
    \end{enumerate}
\end{algorithm}
\begin{proof}
By Lemma \ref{lem:GST_free} every orbit contains a representative with codomain $W_t$, and the subset of such isogenies is a torsor for $\ker(e_T) / (\ker(e_S) \cap \ker(e_T))$.  Correctness now follows from the correctness of Algorithm \ref{alg:min_isog}.
\end{proof}

We use Algorithm~\ref{alg:GST_min} to compute orbit representatives for isogenies of degree $n$.

\begin{algorithm}\label{alg:GST_orbits}$ $

    Input: An isogeny class of ordinary abelian varieties over a finite field $\F_q$ given by its Frobenius order $R$, and an integer $D>1$.

    Output: For weak equivalence classes $s,t$ and $1 < n \mid D$, the finite set $\cR_{s,t}^{=n}$.
       \begin{enumerate}
        \item Run Algorithm~\ref{alg:GST_min} to compute $\cR_{s,t}^{=n,1}$ for $1 < n \mid D$.  Initialize $\cR_{s,t}^{=n} = \cR_{s,t}^{=n,1}$.
        \item Loop over the divisors $n$ of $D$ in increasing order, and
            for each divisor $1 < d_2 < n$ of $n$, loop over triples of weak equivalence classes $s, t, u$ and $(E_2,E_1) \in \cR_{t,u}^{=d_2,1}\times \cR_{s,t}^{=n/d_2}$:
            \begin{enumerate}
                \item For each $E$ in the composition $E_2 \usq E_1$, use Algorithm \ref{alg:equiv_orbit} to determine if its orbit has already been added to $\cR_{s,u}^{=n}$, and add it if not.
            \end{enumerate}
        \item Return the $\cR^{=n}_{s,t}$.
    \end{enumerate}
\end{algorithm}
\begin{proof}
The correctness of Algorithm \ref{alg:equiv_orbit} guarantees that no two returned isogenies lie in the same orbit, and Proposition \ref{prop:GSTcompose} implies that any composition of minimal isogenies will be in the same orbit as one constructed in $E_2 \usq E_1$.  Finally, every isogeny can be expressed as a composition of minimal isogenies as constructed in Algorithm~\ref{alg:GST_min} by Proposition~\ref{prop:isog} .
\end{proof}

\section{Polarizations}\label{sec:pols}
In this section, we compute the polarizations of elements belonging to a given isogeny class.
Throughout, for an isogeny class $\cI_h$, let $g$ be the dimension of the abelian varieties it contains; then the \'etale algebra $K=\Q(\pi)$ satisfies $[K:\Q] = 2g$.
Recall that $K$ is equipped with an involution $x\mapsto \overline{x}$.
We say that an element $x\in K$ is \emph{totally imaginary} if $x=-\overline x$ and \emph{totally real} if $x=\overline{x}$.
A totally real element $x\in K^\times$ is called \emph{totally positive} if $\phi(x)>0$ for every homomorphism $\phi\colon K\to \C$.
Let $\Phi$ be a CM-type of $K$, that is, a choice of $g$ homomorphisms $K\to \C$, one per conjugate pair.
We say that $x\in K$ is \emph{$\Phi$-positive} if $\Im(\phi(x))>0$ for every $\phi\in \Phi$.
For each abelian variety $A$ we denote by $A^\vee$ its dual.  Given a weak equivalence class $s$, we denote by $s^\vee$ the weak equivalence class containing the dual of each element of $s$.

We assume that $\Phi$ satisfies the Shimura--Taniyama formula; see \cite[Notation~4.6]{Howe95} and \cite[2.1.4.1]{chaiconradoort14}.
Let $A,A_1,A_2 \in \cI_h$, and set $I =\cF(A)$, $I_1 = \cF(A_1)$ and $I_2 = \cF(A_2)$.
Let $\varphi \colon A\to A^\vee$ be an isogeny and $(A_1,\psi_1),(A_2,\psi_2)$ be polarized abelian varieties.
Set $\mu = \cF(\varphi)$, $\lambda_1 = \cF(\psi_1)$ and $\lambda_2 = \cF(\psi_2)$.

\begin{proposition}[{\cite[Theorem~5.4]{MarAbVar18}}]\label{prop:pols}
    The following statements hold:
    \begin{enumerate}[(i)]
        \item \label{prop:pols:dual} $\cF(A^\vee)=\overline I^\dagger$.
        \item \label{prop:pols:is_pol} $\varphi$ is a polarization if and only if $\mu$ is totally imaginary and $\Phi$-positive.
        \item \label{prop:pols:is_isom} $(A_1,\psi_1)$ and $(A_2,\psi_2)$ are isomorphic if and only if there is an element $y\in K^\times$ satisfying $yI_1 = I_2$ and $\lambda_1=\overline{y}\lambda_2 y$.
    \end{enumerate}
\end{proposition}

We break the search for polarized abelian varieties up into two steps.  First, we give two different algorithms for producing isogenies of the form $A \to A^\vee$, one based on Algorithm \ref{alg:compositions} and one on Algorithm \ref{alg:GST_orbits}.  Second, we give an algorithm that produces the list of all polarized abelian varieties with polarization of degree dividing $D$ within an isogeny class given the set of such $A \to A^\vee$ as input.

\begin{algorithm}\label{alg:dualiso_from_GD}$ $

Input: An isogeny class $\cI_h$ and an integer $D> 1$. 
Output: Isogenies $\varphi \colon  A \to A^\vee$ of degree dividing $D$, one per equivalence class.

\begin{enumerate}
    \item Use Algorithm \ref{alg:compositions} to compute the isogeny graph $\cG^D$.
        \item For each vertex $V$ do:
            \begin{enumerate}
            \item \label{alg:dualiso_from_GD:Vv} Identify the vertex $V^\vee$ corresponding to the ideal class of $\overline{I_{V}}^\dagger$.
        \item \label{alg:dualiso_from_GD:iota} Compute $\iota \in K^\times$ such that $\iota I_{V^\vee}=\overline{I_{V}}^\dagger$.
        \item For each edge in $\cG^D$ represented by $x_0I_V \subseteq I_{V^\vee}$, yield the isogeny represented by $\iota x_0I_V\subseteq \overline{I_{V}}^\dagger$.
    \end{enumerate}
\end{enumerate}
\end{algorithm}

\begin{proof}
The correctness follows from Proposition \ref{prop:pols}.\ref{prop:pols:dual} and the definition of $\cG^D$.
\end{proof}

\begin{remark}
    In step \ref{alg:dualiso_from_GD:Vv}, one needs to solve a discrete logarithm problem in $\ICM(R)$.
    This expensive step can be avoided by tracking the action of duality when constructing the set of vertices of $\cG$ as orbits of Picard groups; see Proposition~\ref{prop:pic_weak}.
    By Lemma~\ref{lem:dagger}, this action is encoded by the formula $\overline{(IL)}^\dagger = \overline{I}^\dagger\overline{L}^{-1}$.
\end{remark}

Alternatively, we can iterate over representatives from Algorithm \ref{alg:GST_orbits}, which saves time when $\Pic(\ZFV)$ is large.

\begin{algorithm}\label{alg:dualiso_from_Rst} $ $

    Input: An isogeny class $\cI_h$ and an integer $D>1$.

Output:  Isogenies $\varphi \colon  A \to A^\vee$ of degree dividing $D$, one per equivalence class.

\begin{enumerate}
\item Use Algorithm \ref{alg:GST_orbits} to compute $\cR_{s,s^\vee}^n$ for $1 < n \mid D$.
\item For each weak equivalence class $s$, compute $J_s = (\overline W_{s^\vee}^\dagger : W_s)$; set $\frj_s = [J_s]$.
\item For each endomorphism ring $S$, compute the endomorphism $\beta$ of the finite abelian group $G_{S, \overline{S}}$ defined by $\beta(\frb) = \frb \overline{\frb}$.
\item For each weak equivalence class $s$ and each $\varphi \in \cR_{s,s^\vee}^D$ do:
\begin{enumerate}
\item Let $S$ be the multiplicator ring of $s$, $\frw_s \fra$ be the domain of $\varphi$ and $\frw_{s^\vee}$ be the codomain.  Choose an element $\fra' \in \image(\beta)$ mapping to $\fra^{-1} \frj_s$ in $\Pic(S)$ and $\overline{\fra}^{-1}\overline{\frj}_s$ in $\Pic(\overline{S})$; if no $\fra'$ exists then continue to the next $\varphi$.  Otherwise, choose $\frb_0 \in G_{S,\overline{S}}$ with $\beta(\frb_0) = \fra'$.
\item For each $\frc \in \ker(\beta)$, let $\frb_0 \frc \cdot \varphi$ be represented by $x_0I \subseteq J$ and:
\begin{enumerate}
\item Compute $\iota \in K^\times$ such that $\iota J = \overline I^\dagger$, and yield $\iota x_0 I \subseteq \overline{I}^\dagger$.
\end{enumerate}
\end{enumerate}
\end{enumerate}
\end{algorithm}

\begin{proof}
The $G_{S,\overline{S}}$-orbit of $\varphi$ consists of isogenies with domain $\frw_s \fra \frb$ and codomain $\frw_{s^\vee} \frb$ for $\frb \in G_{S,\overline{S}}$.  By Proposition~\ref{prop:pols}.\ref{prop:pols:dual} and Lemma~\ref{lem:dagger}, we seek $\frb$ so that
\begin{align*}
\frw_s \fra \frb &= \overline{(\frw_{s^\vee} \frb)}^\dagger = \overline{\frw_{s^\vee}}^\dagger \overline{\frb}^{-1} = \frw_s \frj_s \overline{\frb}^{-1} \\
\frw_{s^\vee} \overline{\frj}_s \overline{\fra}^{-1} \overline{\frb}^{-1} = \overline{\frw_s}^\dagger \overline{\fra}^{-1} \overline{\frb}^{-1} = \overline{(\frw_s \fra \frb)}^\dagger &= \frw_{s^\vee} \frb
\end{align*}
which will occur when $\frb \overline{\frb} = \fra'$ in $\Pic(S)$ (by the first line) and in $\Pic(\overline{S})$ (by the second line).  The elements $\frb_0 \frc$ are precisely the solutions to this equation.  So far we have worked up to isomorphism; to get an isogeny of the form $\varphi \colon  A \to A^\vee$ we need to compose with $\iota$ so that the codomain is precisely the dual of the domain, not merely isomorphic to it.
\end{proof}

\begin{algorithm}\label{alg:pols_from_dualiso}$ $

Input: An isogeny class $\cI_h$ with Frobenius order $R$, an integer $D>1$, and a CM-type $\Phi$ of $K = \Q[x]/(h(x))$ satisfying the Shimura--Taniyama formula.
Output: Pairs $(I,\mu)$ representing the isomorphism classes of polarized abelian varieties in $\cI_h$ with polarizations of degree dividing $D$.

\begin{enumerate}
    \item\label{alg:pols_from_dualiso:get_isog} Use either Algorithm \ref{alg:dualiso_from_GD} or \ref{alg:dualiso_from_Rst} to compute a list $\cL$ of isogenies $\varphi \colon  A \to A^\vee$ of degree dividing $D$.
    \item For each overorder $S$ of $R$ do:
    \begin{enumerate}
        \item Compute a transversal $\mathcal{S}$ of the quotient $S^\times/S^\times_+$, where $S^\times_+$ is the subgroup of $S^\times$ consisting of totally real and totally positive units.
        \item Compute representatives $v_1,\ldots,v_n$ in $S^\times$ of $\dfrac{S^\times_+}{S^\times_+ \cap \left\langle w\overline{w}:w \in S^\times \right\rangle}$.
    \end{enumerate}
    \item For each isogeny $x I_A \subseteq \overline{I_A}^\dagger$ in $\cL$ do:
    \begin{enumerate}
        \item Set $S=(I_A:I_A)$ and assign $\mathcal{S}$ and $v_1,\ldots,v_n$ as above.
        \item \label{alg:pols_from_dualiso:is_pol} If there exists $v\in \mathcal{S}$ with $\mu = xv$ totally imaginary and $\Phi$-positive then yield the non-isomorphic polarizations $(I_A,\mu v_1),\ldots,(I_A,\mu v_n)$.
    \end{enumerate}
\end{enumerate}
\end{algorithm}

\begin{proof}
    By Proposition~\ref{prop:pols}.\ref{prop:pols:is_pol}, each element $\mu v_i$ represents a polarization of $I_A$.
    Each $\mu v_i$ is an isogeny from $I_A$ to $\overline{I}_{A}^\dagger$ of degree dividing $D$, since it is a composition of automorphisms $v$ and $v_i$ of $I_A$ and an isogeny $x$ of degree dividing $D$.

    We now prove that the output consists of non-isomorphic polarized abelian varieties.
    Pick $(I_{A_1},\mu_1),(I_{A_2},\mu_2)$ in $\mathcal{L}$.
    Since the ideal $I_A$ representing the vertex of $A$ is deterministically chosen,
    we can assume that $A_1=A_2$.
    Then Proposition~\ref{prop:pols}.\ref{prop:pols:is_isom} says that $(I_{A_1},\mu_1)$ and $(I_{A_1},\mu_2)$ are isomorphic if and only if there is a unit $u\in S^\times$ such that $\mu_1=\overline{u}\mu_2u$.
    A necessary condition is that $\mu_1$ and $\mu_2$ come from the same edge $x_0$ by Corollary~\ref{cor:eq_cl_toJ}.\ref{cor:eq_cl_toJ:4} (composed with an isomorphism $\iota$ as in Algorithm \ref{alg:dualiso_from_GD} or \ref{alg:dualiso_from_Rst}).
    There exists at most one element $v$ of $\mathcal{S}$ such that the element $\mu=\iota x$ is totally imaginary and $\Phi$-positive.
    By writing $\mu_1=\mu v_i$ and $\mu_2=\mu v_j$, we see that $\mu_1$ and $\mu_2$ give rise to isomorphic polarized abelian varieties if and only if $\mu_1=\mu_2$.

It remains to show that any $(I,\lambda)$ with $\deg(\lambda) \mid D$ is isomorphic to some $(I_A,\mu v_i)$ in the output.
    Proposition~\ref{prop:pols}.\ref{prop:pols:is_isom} states that this happens precisely if there exists $x\in K^\times$ such that $xI_A=I$ and $\mu v_i=\overline{x}\lambda x$.
    Now, since $\lambda$ is an isogeny of degree dividing $D$, there exists a unique edge $x_0I_A \subseteq I_{A^\vee}$ together with elements $y,z\in K^\times$ such that $yI_A=I$, $z\overline{I}^\dagger=I_{A^\vee}$ and $x_0=z\lambda y$ by Proposition~\ref{prop:same_edge}.\ref{prop:same_edge:2}.
    Take $\iota \in K^\times$ satisfying $\iota I_{A^\vee} = \overline{I_A}^\dagger$.
    Let $\overline{u}=\overline{y}/z\iota$.
    Since $\overline{y}\overline{I}^\dagger = z\iota\overline{I}^\dagger$, then $u\in S^\times$.
    We summarize the situation in the following diagram:
    \[\begin{tikzcd} [sep = .5 cm]
        I\dar["\lambda"] & I_A \lar["y"]\dar["x_0"] & \\
        \overline{I}^\dagger\drar["\overline{y}"]\rar["z"] & I_{A^\vee}\rar["\iota"] & \overline{I_A}^\dagger \dlar["\overline{u}"]\\
         & \overline{I_A}^\dagger &
    \end{tikzcd}\]
    We have $\overline{u}^{-1}\overline{y}\lambda y = \iota x_0$.
    By multiplying by $u^{-1}$, we see that $\iota x_0 u^{-1}$ is totally imaginary and $\Phi$-positive, since $\lambda$ is so.
    Let $v$ be the representative in $\mathcal{S}$ of $u^{-1}$, so that $vu\in S^\times_+$.
    Let $v_i$ be the representative of $(vu)^{-1}$ in ${S^\times_+}/({S^\times \cap \left\langle w\overline{w}:w \in S^\times \right\rangle})$.
    Hence there exists $w\in S^\times$ such that $vuv_i=w\overline{w}$.
    By setting $x=wu^{-1}y$, we get $\overline{x}\lambda x=\mu v_i$, that is, an isomorphism between $(I_A,\mu v_i)$ and $(I,\lambda)$.
\end{proof}

\subsection*{Acknowledgements}\
Costa was supported by Simons Foundation grant SFI-MPS-Infrastructure-00008651, Roe by SFI-MPS-SSRFA-00005427, and both by 550033.
Dupuy was supported by National Science Foundation grant DMS-2401570.
Dupuy and Vincent are grateful for the hospitality of the GAATI laboratory at the Universit\'{e} de Polyn\'{e}sie fran\c{c}aise, where they completed this work.
Marseglia was supported by NWO grant VI.Veni.202.107,
and by Marie Sk{\l}odowska-Curie Actions - Postdoctoral Fellowships 2023 (project 101149209 - AbVarFq).
Vincent was supported by a Simons Foundation Travel Support for Mathematicians grant.

\printbibliography

\end{document}